\documentclass[10pt]{article}

\usepackage{amssymb}
\usepackage[font={small, it}]{caption}
\usepackage{amsmath}
\usepackage{floatrow}
\usepackage{halloweenmath}
\usepackage{times}
\usepackage{ stmaryrd }

\usepackage{amsthm}
\usepackage{authoraftertitle}
\usepackage[
backend=biber,
style=numeric,
sorting=nty
]{biblatex}
\usepackage{relsize}
\usepackage{xcolor}
\usepackage{mathrsfs}
\usepackage[colorlinks = true,
            linkcolor = black,
            urlcolor  = blue,
            citecolor = black,
            hypertexnames = false,
            anchorcolor = blue]{hyperref}
\usepackage{wasysym}
\usepackage[letterpaper, portrait, margin=1in]{geometry}
\usepackage{graphicx}
\usepackage{tikz}
\usepackage{tikz-3dplot}
\usepackage{pgfplots}
\usetikzlibrary{decorations.pathmorphing,patterns}
\usetikzlibrary{decorations.markings}
\usepackage{lipsum}
\usepackage{float}
\usepackage{subcaption}
\usepackage[object=vectorian]{pgfornament}
\usepackage{mwe}
\usepackage{bigints}
\usepackage{titlesec}
\titleformat{\paragraph}
{\normalfont\normalsize\bfseries}{\theparagraph}{1em}{}
\titlespacing*{\paragraph}
{0pt}{3.25ex plus 1ex minus .2ex}{1.5ex plus .2ex}
\usepackage{mathtools}
\usepackage{pgfplots}
\pgfplotsset{compat=1.15}
\usepackage{lastpage}
\usepackage{enumitem}
\usepackage{tensor}
\usepackage{mathtools}
\usepackage{fancyhdr} 
\usepackage{setspace}
\usepackage{tikz}
\usetikzlibrary{hobby}

\usepackage{pst-node}

\usepackage{tikz-cd}
\usepackage{thm-restate}
\usepackage{thmtools}

\usepackage{cleveref}

\usepackage{apptools}
\AtAppendix{\counterwithin{lemma}{section}}

\usepackage{setspace}

\theoremstyle{definition}

\theoremstyle{plain}
\declaretheorem[name=Theorem]{theorem}

\newtheorem{lemma}[theorem]{Lemma}

\newtheorem{remark}{Remark}

\newtheorem{proposition}[theorem]{Proposition}

    \newcommand{\ve}{\varepsilon}
    \newcommand{\Ocal}{\mathcal{O}}
    
    \newcommand{\al}{\alpha}
    \newcommand{\be}{\beta}

    \newcommand{\lahat}{\hat{\la}}

    \newcommand{\Gbar}{\,\overline{\!G}}
    
    \newcommand{\De}{\Delta}

    \newcommand{\Fbar}{\,\overline{\!F}}
    
    \newcommand{\Dbar}{\,\overline{\!D}}

    \newcommand{\Kbar}{\,\overline{\!K}}

    \newcommand{\an}{\mathrm{an}}

    \newcommand{\albar}{\bar{\al}}

    \newcommand{\Pbb}{\mathbb{P}}
    
    \newcommand{\fbar}{\bar{f}}
    \newcommand{\gbar}{\bar{g}}

    \newcommand{\de}{\delta}

    \newcommand{\xbar}{\bar{x}}
    \newcommand{\abar}{\bar{a}}
    \newcommand{\ybar}{\bar{y}}

    \newcommand{\ze}{\zeta}

    \DeclareMathOperator{\Disc}{Disc}

    \newcommand{\vp}{\varphi}

    \DeclareMathOperator{\chara}{char}
    
    \newcommand{\Pfrak}{\mathfrak{P}}
    \newcommand{\mfrak}{\mathfrak{m}}

    \newcommand{\wvec}{\mathbf{w}}

    \newcommand{\vvec}{\mathbf{v}}

    \newcommand{\Zbb}{\mathbb{Z}}

    \newcommand{\bebar}{\bar{\be}}

    \newcommand{\Qbb}{\mathbb{Q}}
    \newcommand{\ga}{\gamma}
    \newcommand{\la}{\lambda}
    \newcommand{\RomanNumeralCaps}[1]
        {\MakeUppercase{\romannumeral #1}} 
    
    \newcommand{\Rbb}{\mathbb{R}}

    \newcommand{\Acal}{\mathcal{A}}
    \newcommand{\pfrak}{\mathfrak{p}}

    \newcommand{\Cbb}{\mathbb{C}}

\makeatletter
\DeclareRobustCommand{\cev}[1]{%
  {\mathpalette\do@cev{#1}}%
}
\newcommand{\do@cev}[2]{%
  \vbox{\offinterlineskip
    \sbox\z@{$\m@th#1 x$}%
    \ialign{##\cr
      \hidewidth\reflectbox{$\m@th#1\vec{}\mkern4mu$}\hidewidth\cr
      \noalign{\kern-\ht\z@}
      $\m@th#1#2$\cr
    }%
  }%
}
\makeatother

\makeatletter
\DeclarePairedDelimiterX{\pmodx}[1]{(}{)}{{\operator@font mod}\mkern6mu#1}
\renewcommand{\pmod}{%
  \allowbreak
  \if@display\mkern18mu\else\mkern8mu\fi
  \pmodx
}
\makeatother

\DeclarePairedDelimiterX\braket[2]{\langle}{\rangle}{#1 \delimsize\vert #2}

\makeatletter
\newcommand{\colim@}[2]{%
  \vtop{\m@th\ialign{##\cr
    \hfil$#1\operator@font colim$\hfil\cr
    \noalign{\nointerlineskip\kern1.5\ex@}#2\cr
    \noalign{\nointerlineskip\kern-\ex@}\cr}}%
}
\newcommand{\colim}{%
  \mathop{\mathpalette\colim@{\rightarrowfill@\scriptscriptstyle}}\nmlimits@
}
\newcommand{\dlim}{%
  \mathop{\mathpalette\varlim@{\rightarrowfill@\scriptscriptstyle}}\nmlimits@
}
\newcommand{\inlim}{%
  \mathop{\mathpalette\varlim@{\leftarrowfill@\scriptscriptstyle}}\nmlimits@
}

\tikzset{middlearrow/.style={decoration={markings, mark= at position 0.52 with {\arrow{#1}},},postaction={decorate}}}

\usepackage{authblk}

\title{Infinitely wildly ramified arboreal representations for postcritically finite polynomials with potential good reduction}
\author[a]{Alex Feiner\footnote{Corresponding author. Email address: alexander\_feiner@brown.edu. ORCID iD: 0000-0001-9588-7860.}}
\affil[a]{Department of Mathematics, Brown University, 151 Thayer St, Providence, RI 02912.}
\date{}
\begin{document}
\maketitle

\begin{abstract}
Let $K$ be a number field, let $v$ be a finite place of $K$, let $f\in K[z]$ be a degree $d\geqslant2$ polynomial with $v|d$, and let $a\in K$. We show that if $f$ is postcritically bounded and has potential good reduction with respect to $v$, then the arboreal representation associated to the pair $(f,a)$ is either finite or infinitely wildly ramified above $v$.

\end{abstract}

\section{Introduction}
\hspace{\parindent} Let $K$ be a number field. For $f\in K[z]$ with degree $d\geqslant2$, let $f^n$ denote the $n$th iterate of $f$ (with the convention that $f^0(z)=z$), and for $a\in K$, let $f^{-n}(a)$ denote the set of solutions in $\Kbar$ to $f^n(z)=a$. We define the set of iterated preimages of $a$ under $f$ to be
    \[T_\infty(f,a)=\{a\}\cup f^{-1}(a)\cup f^{-2}(a)\cup\cdots.\]
Here, we are considering $T_\infty(f,a)$ as a subset of $\Kbar$, and not as a tree or a multiset. Note that if $a$ is not in the postcritical set of $f$ then $f^n(z)-a$ has $d^n$ distinct roots in $\Kbar$ for all $n\geqslant0$, which shows that $T_\infty(f,a)$ is infinite. If $v$ is a finite place of $K$, then we say that $f$ has \emph{good reduction} at $v$ if $f$ is monic and the valuation of all the coefficients of $f$ are non-negative. We say that $f$ has \emph{potential good reduction} at $v$ if there is some degree $1$ polynomial $\mu\in\Kbar[z]$ for which $\mu\circ f\circ\mu^{-1}\in\Kbar[z]$ has good reduction with respect to any extension of $v$ to $\Kbar$. We say that $f$ is \emph{postcritically bounded} with respect to $v$ if for any extension of $v$ to $\Kbar$, the set 
\[
\{f^n(\ga) \mid n\geqslant0, \ga \text{ a critical point of }f\}
\quad\text{is $v$-adically bounded.}
\]
\par
Define field extensions $K_n/K$ and $K_\infty/K$ by setting 
\[
K_n=K(f^{-n}(a))
\quad\text{and}\quad
K_\infty=K(T_\infty(f,a))
=\bigcup_{n\ge0} K_n.
\]
We say that $K_\infty/K$ is \emph{infinitely ramified} above a prime $\pfrak$ of $K$ if for every sequence of primes 
    \[\pfrak\subseteq\Pfrak_1\subseteq\Pfrak_2\subseteq\cdots,\]
where $\Pfrak_n$ is a prime of $K_n$, we have that $e(\Pfrak_n|\pfrak)\to\infty$ as $n\to\infty$, and we say that $K_\infty/K$ is \emph{infinitely wildly ramified} above $\pfrak$ if $v_p(e(\Pfrak_n|\pfrak))\to\infty$ as $n\to\infty$, where $p\in\Zbb$ is the unique rational prime lying below $p$. Note that the fact that each extension $K_n/K$ is Galois shows that if this is true for one such sequence of primes then it is true for any such sequence of primes. In \cite[Conjecture 6]{FinRamiForPreim}, they conjecture that if $f$ is a PCF polynomial of degree $d\geqslant2$ and the extension $K_\infty/K$ is an infinite extension, then there is at least one prime $\pfrak$ of $K$ lying above $d$ for which $K_\infty/K$ is infinitely wildly ramified above $\pfrak$ (see also \cite[conjecture 1.2]{TamePropGalGroups} and \cite[section 1]{ImageArborealGalReps}). We prove that this is actually true in the case where $f$ has potential good reduction at a finite place $v\mid d$ of $K$. In particular, we show the following:

\begin{restatable}[]{theorem}{InfinitelyRamified}
\label{thm:InfinitelyRamified}
Let $K$ be a number field and let $v$ be a finite place of $K$ lying above the rational prime $p\in\Zbb$. Suppose that $f\in K[z]$ is a degree $d$ polynomial that is postcritically bounded at $v$ and has potential good reduction at $v$, and that $p|d$. Let $a\in K$ be such that $\#T_\infty(f,a)=\infty$, let $K_n=K(f^{-n}(a))$, and let $K_\infty$ be the union of the $K_n$ for $n\geqslant0$. Then $K_\infty/K$ is infinitely wildly ramified above $\pfrak$, where $\pfrak$ is the prime of $K$ corresponding to $v$. 
\end{restatable}

We will see that we can pass to the completion $K_v$ in order to translate the above into a problem about local fields of characteristic $0$ and residue field characteristic $p$, where $v\mid p$. In the case where $K$ is a local field of characteristic $0$ and residue field characteristic $p$, if we have a tower of finite Galois extensions $\cdots /K_2/K_1/K$ and $K_\infty$ is the union of all the $K_n$, then we say that $K_\infty/K$ is \emph{infinitely ramified} if $e(K_n/K)\to\infty$ as $n\to\infty$, and we say that $K_\infty/K$ is \emph{infinitely wildly ramified} if $v_p(e(K_n/K))\to\infty$ as $n\to\infty$. We will show that Theorem \ref{thm:InfinitelyRamified} is implied by the following:

\begin{restatable}[]{theorem}{InfinitelyRamifiedLocalField}
\label{thm:InfinitelyRamifiedLocalField}
Let $(K,v)$ be a local field of characteristic $0$ and residue field characteristic $p$. Let $f\in\Ocal_K[z]$ be a monic postcritically bounded polynomial that fixes $0$ and has degree $d$, with \text{$p\mid d$}. Suppose that $a\in\Ocal_K$ is not in the postcritical set of $f$. Let $K_n=K(f^{-n}(a))$ and let $K_\infty$ denote the union of all the $K_n$ for $n\geqslant0$. Then $K_\infty/K$ is infinitely wildly ramified. 
\end{restatable}

Often, a postcritically bounded polynomial is forced to have potential good reduction at some finite place $v|d$ for geometric reasons. For example, if $f\in K[z]$ is a postcritically bounded polynomial of prime power degree $d=p^\ell\geqslant2$, then it has potential good reduction at $v$ (see Lemma \ref{lem:wildplacegoodreduction'-prime-power-degree}). Thus Theorem 1 shows that if $f\in K[z]$ is a polynomial of prime power degree $d=p^\ell\geqslant2$ that is postcritically bounded at $v$, and if $a\in K$ is such that $\#T_\infty(f,a)=\infty$, then $K_\infty/K$ is infinitely wildly ramified above $\pfrak$. In fact, we will see that an appropriate conjugate of $f$ will have monomial good reduction, which allows us to recover a special case of a result due to Sweeney (see \cite[Theorem 4.3]{SweeneySensThm}) concerning arboreal representations of prime power degree polynomials. However, work by Anderson, Manes, and Tobin (see \cite[Theorem 1]{BelyiPolys}) shows that it is possible for a postcritically finite polynomial to fail to have potential good reduction at some $v|d$ in the case where $d$ is not a prime power. Thus the issue of whether the results of Theorems \ref{thm:InfinitelyRamified} and \ref{thm:InfinitelyRamifiedLocalField} still hold in the case where $f$ does not have potential good reduction at $v$ remains a natural question emerging from this paper. Finally, we note that the problem of wild ramification in arboreal representations has been studied for the case of unicritical polynomials in both \cite{LocalArborealReps} and \cite{familyoflowdegreeextensions}.

\subsection{Outline}
\hspace{\parindent} We begin by examining the behavior of images and preimages of disks in $\Cbb_v$ under a monic polynomial $f$ with integral coefficients. In particular, we show that the image of large enough disks contained in the ring of integers of $\Cbb_v$ behave very nicely, except when they are centered at elements contained in certain ``bad directions" at the Gauss point. We then use an explicit formula for the valuation of the discriminant of iterates of $f$, as well as a result proved in the appendix, to create a lower bound on the ramification indices of the field extensions $K(f^{-n}(a))/K$. Using this, we relate the results about images of disks in $\Cbb_v$ to the ramification of the iterated field extensions in question, and show that $K_\infty/K$ is infinitely wildly ramified if $T_\infty(f,a)$ does not contain elements lying in these ``bad directions." We then reduce the iterated set of preimages of $a$ under $f$ modulo the maximal ideal of the ring of integers of $\Cbb_v$, and analyze this in order to show that $K_\infty/K$ is infinitely wildly ramified even when $T_\infty(f,a)$ includes elements lying in these ``bad directions."

\subsection{Notation}
\hspace{\parindent} Throughout this paper, unless explicitly stated otherwise, we use the following notation:
\begin{enumerate}[label=(\arabic*)]
    \item $(K,v)$ is a local field of characteristic $0$ and residue field characteristic $p$. Furthermore, unless explicitly stated otherwise, we assume that $v=v_\pfrak$ for the prime ideal $\pfrak$ of $K$. 
    \item $\Cbb_v$ is the completion of an algebraic closure of $K$. 
    \item For $a\in\Cbb_v$ and $r\geqslant0$, we let $D(a,r),\Dbar(a,r)\subseteq\Cbb_v$ denote the open and closed disks of radius $r$ centered at $a$ in $\Cbb_v$, respectively, given by
        \[D(a,r)=\{z\in\Cbb_v \mid |z-a|_v<r\},\qquad \Dbar(a,r)=\{z\in\Cbb_v \mid |z-a|_v\leqslant r\}.\]
    We use the words ``open" and ``closed" solely to refer to how the disk $D$ is defined, and not to refer to any topological properties of the disk. 
    \item For $a\in\Cbb_v$ and $0<r_1<r_2$, we let $\Acal(a; r_1,r_2)\subseteq\Cbb_v$ denote the open annulus centered at $a$ in $\Cbb_v$, given by
        \[\Acal(a; r_1,r_2)=\{z\in\Cbb_v \mid r_1<|z-a|_v<r_2\}.\]
    \item We let $\Ocal=\Dbar(0,1)$ denote the ring of integers of $\Cbb_v$, $\mfrak=D(0,1)$ be its maximal ideal, and $k=\Ocal/\mfrak$ be its residue field with $\chara (k)=p$.    
    \item If $a\in\Ocal$ then we denote by $\abar\in k$ its reduction modulo $\mfrak$, and similarly for $f\in\Ocal[z]$. 
\end{enumerate}

\section*{Acknowledgements}
\hspace{\parindent} The author would like to thank Nicole Looper for countless extremely helpful conversations surrounding the material. They would also like to thank Joe Silverman for valuable conversations about the material and advice in editing the paper, and Rob Benedetto for very useful correspondences about the material.

\section{Non-Archimedean Dynamics}
\hspace{\parindent} We let $\Pbb^1_\an$ denote the Berkovich projective line associated to $\Cbb_v$ (see \cite[Chapter 6]{RobBook}). For $a\in\Cbb_v$ and $r>0$, we denote the type \RomanNumeralCaps{2} or \RomanNumeralCaps{3} point in $\Pbb^1_\an$ corresponding to the closed disc $\Dbar(a,r)$ by $\ze(a,r)$. For $\ze=\ze(a,r)\in\Pbb^1_\an$ such a point and $x\in\Pbb^1_\an\setminus\{\ze\}$, we let $\vvec_\ze(x)\subseteq\Pbb^1_\an\setminus\{\ze\}$, called the direction at $\ze$ containing $x$, denote the connected component of $\Pbb^1_\an\setminus\{\ze\}$ that contains $x$. In particular, if $\ze=\ze(0,1)$ is the Gauss point and $x\in\Ocal$, then 
    \[\vvec_{\ze(0,1)}(x)\cap\Cbb_v=D(x,1).\]
For $f\in\Cbb_v[z]$, $\ze=\ze(a,r)\in\Pbb^1_\an$, and $\vvec=\vvec_\ze(a)$ a direction at $\ze$ not containing $\infty$, we let $\deg_{\ze,\vvec}(f)$ denote the degree of $f$ at $\ze$ in the direction of $\vvec$. That is, there is some $0<\la<1$ for which $f$ has no roots on the annulus $\Acal(a; \la r,r)\subseteq\Cbb_v$, and we let $\deg_{\ze,\vvec}(f)$ denote the common inner and outer Weierstrass degree of $f$ on this annulus. Furthermore, we let $f_\sharp(\vvec)$ denote the direction at $f(\ze)$ that contains $f(\Acal(a;\la r,r))$ (see \cite[definition 7.18]{RobBook}).

We now state two propositions about preimages and images of disks in $\Cbb_v$ under polynomials that will be crucial throughout this document. In general, the preimage of a disk in $\Cbb_v$ under a polynomial will behave very nicely. In fact, we have the following:

\begin{lemma}[{\cite[Theorem 3.23]{RobBook}}]\label{lem:Robbook-theorem-3.23}
Let $D\subseteq\Cbb_v$ be a disk and let $f\in\Cbb_v[z]$ be a degree $d$ polynomial. Then $f^{-1}(D)$ is a disjoint union of disks,
    \[f^{-1}(D)=D_1\cup\cdots\cup D_n,\]
where each $D_i$ is open if $D$ is open and closed if $D$ is closed. Furthermore, for each $i=1,\dots,n$ there is an integer $1\leqslant d_i\leqslant d$ such that $f$ maps $D_i$ $d_i$-to-$1$ (counting multiplicity) onto $D$, and $d_1+\cdots+d_n=d$. 
\end{lemma}

In the notation of Lemma \ref{lem:Robbook-theorem-3.23}, we will refer to the $D_i$ as the \emph{disk components} of $f^{-1}(D)$. It is easy to see that if $f(D_0)=D$ then $D_0$ is one of the disk components of $f^{-1}(D)$. There are also nice results about how a polynomial $f\in\Cbb_v[z]$ stretches the distance between points in a disk in $\Cbb_v$.

\begin{lemma}[{\cite[proposition 3.20]{RobBook}}]\label{lem:Robbook-proposition-3.20}
Let $D\subseteq\Cbb_v$ be a disk of radius $r$, let $f\in\Cbb_v[z]$, and suppose that $f(D)$ is a disk of radius $s$. Then for all $x,y\in D$,
    \[|f(x)-f(y)|_v\leqslant\frac{s}{r}|x-y|_v.\]
\end{lemma}

In particular, if we take $f\in\Ocal[z]$, so $f(\Ocal)\subseteq\Ocal$, then Lemma \ref{lem:Robbook-proposition-3.20} shows that
    \[|f(x)-f(y)|_v\leqslant|x-y|_v\]
for all $x,y\in\Ocal$. Another immediate consequence of Lemma \ref{lem:Robbook-proposition-3.20} is the following:
\begin{lemma}\label{lem:preimage-of boundary-is-on-boundary}
Suppose that $f\in\Cbb_v[z]$, that $D$ is a closed disk of radius $r>0$, and that $f(D)$ is a closed disk of radius $s$. If $\al,\be\in D$ and $|f(\al)-f(\be)|_v=s$, then $|\al-\be|_v=r$. 
\end{lemma}
\begin{proof}
We have that
    \[s=|f(\al)-f(\be)|_v\leqslant\frac{s}{r}|\al-\be|_v,\]
so $|\al-\be|_v\geqslant r$, and thus $|\al-\be|_v=r$. 
\end{proof}

Using these facts about how images of disks in $\Cbb_v$ behave under polynomials, we can show that the image of large enough disks under certain monic polynomials take a particularly nice form. 

\begin{lemma}\label{lem:image-of-disc-under-monic-poly}
Let $f\in\Ocal[z]$ be a degree $d$ monic polynomial. Let $x\in\Ocal$, and suppose that $0<\la<1$ is such that $f-f(x)$ has no roots on $\Acal(x; \la, 1)$. Then $f-f(x)$ has common inner and outer Weierstrass degree $d_x$ on $\Acal(x;\la,1)$, and for any $\la<s<1$, 
    \[f\!\left(\Dbar(x,s)\right)=\Dbar\!\left(f(x),s^{d_x}\right).\] 
\end{lemma}
\begin{proof}
We first make the following observation: For $y\in\Cbb_v$, let $a_{n,y}\in\Cbb_v$ be such that 
    \[f(z)=\sum_{n=0}^{d}a_{n,y}(z-y)^n,\]
so $a_{0,y}=f(y)$ and $a_{d,y}=1$. Then for any $y,y'\in\Cbb_v$, 
\begin{align*}
    f(z)&=\sum_{n=0}^{d}a_{n,y'}(z-y')^n\\
    &=\sum_{n=0}^{d}a_{n,y'}(z-y+(y-y'))^n\\
    &=\sum_{n=0}^{d}\sum_{m=0}^{n}a_{n,y'}\binom{n}{m}(z-y)^m(y-y')^{n-m}\\
    &=\sum_{m=0}^{d}\sum_{n=m}^{d}a_{n,y'}\binom{n}{m}(y-y')^{n-m}(z-y)^m\\
    &=\sum_{m=0}^{d}a_{m,y}(z-y)^m,
\end{align*}
so 
    \[a_{m,y}=\sum_{n=m}^{d}\binom{n}{m}a_{n,y'}(y-y')^{n-m}.\]
Setting $y'=0$ then shows that 
    \[a_{m,y}=\sum_{n=m}^{d}\binom{n}{m}a_{n,0}y^{n-m},\]
so if $f\in\Ocal[z]$ and $y\in\Ocal$ then $a_{m,y}\in\Ocal$ for all $m$.

Now, the fact that $f-f(x)$ has no roots on $\Acal(x; \la,1)$ means that its inner and outer Weierstrass degree on this annulus will be the same integer $d_x=1,\dots,d$ by \cite[proposition 3.32]{RobBook}. Hence 
    \[|a_{d_x,x}|_v=\max_{n=1,\dots,d}|a_{n,x}|_v=1,\]
since $f$ is monic by definition and $|a_{n,x}|_v\leqslant1$ since $a_{n,x}\in\Ocal$. Then \cite[Theorem 3.33]{RobBook} shows that if $y\in\Acal(x;\la,1)$ then 
    \[|f(y)-f(x)|_v=\left|a_{d_x,x}(y-x)^{d_x}\right|_v=|y-x|_v^{d_x}.\]
If $\la<s<1$, then \cite[corollary 3.18]{RobBook} shows that
    \[f\!\left(\Dbar(x,s)\right)=\Dbar(f(x),r),\]
where 
    \[r=\sup_{y\in\Dbar(x,s)}|f(y)-f(x)|_v.\]
By \cite[propositon 3.19]{RobBook}, if $(y_n)_{n\geqslant1}$ is a sequence of points in $\Acal(x; \la,s)\subseteq\Acal(x; \la,1)$ for which $|y_n-x|_v\to s$ as $n\to\infty$, then 
    \[r=\sup_{y\in\Dbar(x,s)}|f(y)-f(x)|_v=\lim_{n\to\infty}|f(y_n)-f(x)|_v=\lim_{n\to\infty}|y_n-x|_v^{d_x}=s^{d_x},\]
and thus 
    \[f\!\left(\Dbar(x,s)\right)=\Dbar\!\left(f(x),s^{d_x}\right).\qedhere\]
\end{proof}

As we will see later in Proposition \ref{prop:general-degree-preimage-of-disks}, we can also tell what the preimage of large enough disks contained in $\Ocal$ look like under $f$ if we know what the reduction of $f$ in $k$ looks like.

\section{Valuative Data For Polynomials Postcritically Bounded at $v|d$}
\hspace{\parindent} In the process of providing a lower bound for certain ramification indices, we will need to have an explicit formula for the valuation of the discriminant of iterates of a polynomial. Note that this is similar to the results of Theorem 1.2 in \cite{FinitelyRamifiedIteratedExtensions}.

\begin{lemma}\label{lem:valuation-of-discriminant}
Suppose that $f,g\in K[z]$ are monic polynomials with degrees $d$ and $e$, respectively. Let $\ga_1,\dots,\ga_{d-1}\in\Kbar$ denote the critical points of $f$, and for $n\geqslant1$, define $N_{g,n}\in\Rbb$ by 
    \[N_{g,n}=\sum_{i=1}^{d-1}v(g(f^n(\ga_i))).\]
Let $\De_{g,n}=\Disc(g\circ f^n)$ for $n\geqslant0$. Then
\begin{equation}
    v(\De_{g,n})=d^nv(\Disc(g))+end^nv(d)+\sum_{i=0}^{n-1}d^iN_{g,n-i}.\label{eq:valuation-of-discriminant}
\end{equation}
In particular, if $f,g\in\Ocal_K[z]$, $v(d)>0$, and $v(\ga_i)\geqslant0$ for all $i$, then $v(\De_{g,n})>0$ for all $n$.
\end{lemma}
\begin{proof}
By \cite[Lemma 2.6]{DensityOfPrimes}, the quantities $\De_{g,n}$ satisfy the recurrence relation
    \[\De_{g,n}=\pm\De_{g,n-1}^{d}d^{ed^n}\prod_{i=1}^{d-1}g(f^n(\ga_i))\]
for all $n\geqslant1$. Hence
    \[v(\De_{g,n})=dv(\De_{g,n-1})+ed^nv(d)+N_{g,n},\]
which shows that equation (\ref{eq:valuation-of-discriminant}) is true when $n=1$ (and it is clearly true when $n=0$). If equation (\ref{eq:valuation-of-discriminant}) is true for some $n\geqslant1$, then
\begin{align*}
v(\De_{g,n+1})&=dv(\De_{g,n})+ed^{n+1}v(d)+N_{g,n+1}\\
    &=d\left(d^nv(\Disc(g))+end^nv(d)+\sum_{i=0}^{n-1}d^iN_{g,n-i}\right)+ed^{n+1}v(d)+N_{g,n+1}\\
    &=d^{n+1}v(\Disc(g))+e(n+1)d^{n+1}v(d)+\sum_{i=0}^{n}d^iN_{g,n+1-i},
\end{align*}
so it is true for all $n$. If $f,g\in\Ocal_K[z]$ and $v(\ga_i)\geqslant0$ for all $i$, then $v(g(f^n(\ga_i)))\geqslant0$ for all positive integers $n$ and $i=1,\dots,d-1$. Furthermore, the fact that $g$ is monic shows that all of its roots are in $\Ocal_K$, so $v(\Disc(g))\geqslant0$, which means that if $v(d)>0$ then $v(\De_{g,n})>0$. 
\end{proof}

We now show that if $f\in\Cbb_v[z]$ is a monic polynomial that fixes $0$ and has integral coefficients and is postcritically bounded with good reduction, then the reduction of $f$ in $k$ actually contains no powers of $z$ that are not divisible by $p$.

\begin{lemma}\label{lem:wildplacegoodreduction'-composite-degree} 
Let $f(z)=z^d+a_{d-1}z^{d-1}+\dots+a_1z\in \Ocal[z]$ be a monic postcritically bounded polynomial of degree $d\geqslant2$, and suppose that \text{$p\mid d$}. Then $v(a_i)>0$ for all $i$ such that $v(i)<v(d)$.  
\end{lemma}
\begin{proof}
Suppose that this is not the case, so there is some $i=1,\dots,d-1$ for which $v(i)<v(d)$ and $v(a_i)=0$. Then 
    \[0=-\frac{v(a_i)}{d-i}=\max_{j=1,\dots,d-1}-\frac{v(a_j)}{d-j}.\]
The Newton polygon $N(f')$ of $f'$ will be the lower convex hull of the points $(j-1,v(a_j)+v(j))$, for $j=1,\dots,d-1$, and $(d-1,v(d))$. Hence $N(f')$ has a line segment of slope
    \[\max_{j=1,\dots,d-1}\frac{v(d)-v(a_j)-v(j)}{d-j}\geqslant\frac{v(d)-v(a_i)-v(i)}{d-i}=\frac{v(d)-v(i)}{d-i}\]
(this is where we are using the assumption that $f$ fixes $0$, since we are assuming that $i\neq0$), which means that $f'$ has a root $\ga\in\Kbar$ with
    \[v(\ga)\leqslant\frac{v(i)- v(d)}{d-i}<0.\]
Hence $|\ga|_v>1$, and 
    \[1=|a_i|_v^{1/(d-i)}=\max_{j=1,\dots,d-1}|a_j|_v^{1/(d-j)},\]
so it follows by \cite[Lemma 2.1]{IngramPCF} (this specific lemma is stated only for absolute values on $\Qbb$, but it works equally well in this case) that
    \[\lahat_{f,v}(\ga)=\log|\ga|_v>0.\]
This contradicts the fact that $f$ is postcritically bounded, so we must have $v(a_i)>0$ for all $i$ such that $v(i)<v(d)$. 
\end{proof}

\begin{remark}\label{critical-points-in-unit-disk}
Note that Lemma \ref{lem:wildplacegoodreduction'-composite-degree} also shows that $|\ga|_v\leqslant1$ for all critical points $\ga$ of $f$, for if $|\ga|_v>1$ then $f$ would not be postcritically bounded. Additionally, Lemma \ref{lem:wildplacegoodreduction'-composite-degree} shows that the reduction of $f$ in $k$ is in the form 
    \[\fbar(z)=\Fbar\!\left(z^{p^\ell}\right),\]
where $\ell=v_p(d)$ and $\Fbar\in k[z]$ is non-constant. In particular, if $d=p^\ell$ then this shows that 
    \[\fbar(z)=z^{p^\ell}.\]
\end{remark}

Slightly modifying Lemma \ref{lem:wildplacegoodreduction'-composite-degree} shows the stronger result that if $f\in K[z]$ is any monic postcritically bounded polynomial of prime power degree $d=p^\ell\geqslant2$ that fixes $0$, then the coefficients of $f$ all have strictly positive valuation. 

\begin{lemma}\label{lem:wildplacegoodreduction'-prime-power-degree}
Let $f(z)=z^d+a_{d-1}z^{d-1}+\cdots+a_1z\in K[z]$ be a monic postcritically bounded polynomial of degree $d=p^\ell\geqslant2$. Then $v(a_i)>0$ for all $i=1,\dots,d$. 
\end{lemma}
\begin{proof}
The fact that $d=p^\ell$ shows that $v(i)<v(d)$ for all $i=1,\dots,d-1$. Suppose that there is some $i=1,\dots,d-1$ for which $v(a_i)\leqslant0$, and fix one such value of $i$ that maximizes $-v(a_i)/(d-i)$. Then
    \[0\leqslant-\frac{v(a_i)}{d-i}=\max_{j=1,\dots,d-1}-\frac{v(a_j)}{d-j},\]
so the Newton polygon of $f'$ has a line segment of slope 
    \[\max_{j=1,\dots,d-1}\frac{v(d)-v(a_j)-v(j)}{d-j}\geqslant\frac{v(d)-v(a_i)-v(i)}{d-i}>-\frac{v(a_i)}{d-i}\geqslant0.\]
Hence there is a root $\ga\in\Kbar$ of $f'$ with 
    \[v(\ga)<\frac{v(a_i)}{d-i}\leqslant0,\]
which implies that $|\ga|_v>1$ and
    \[|\ga|_v>|a_i|_v^{1/(d-i)}=\max_{j=1,\dots,d-1}|a_j|_v^{1/(d-j)}.\]
It then follows by \cite[Lemma 2.1]{IngramPCF} that 
    \[\lahat_{f,v}(\ga)=\log|\ga|_v>0,\]
which contradicts the assumption that $f$ is postcritically bounded. Thus $v(a_i)>0$ for all $i=1,\dots,d-1$.  
\end{proof}

\section{Images of Disks Under Postcritically Bounded Polynomials}

\hspace{\parindent} Let $f\in\Ocal[z]$ be a monic polynomial of degree $d$, and let $x\in\Ocal$. Combining Lemmas \ref{lem:Robbook-theorem-3.23} and \ref{lem:image-of-disc-under-monic-poly}, we see that if $s<1$ is sufficiently close to $1$, then we can write $f^{-1}\!\left(\Dbar(x,s)\right)$ as the (not necessarily disjoint) union
\begin{equation}
    f^{-1}\!\left(\Dbar(x,s)\right)=\bigcup_{i=1}^{N}\Dbar\!\left(x_i,s^{1/d_i}\right),\label{eq:preimage-of-disk-bad-union}
\end{equation}
where $f^{-1}(x)=\{x_1,\dots,x_N\}$ and $d_i=d_{x_i}$. If we remove enough disks from the union in (\ref{eq:preimage-of-disk-bad-union}) then it will become disjoint, but we will still know nothing about what the possible values of the $d_i$ still present in the union are. However, if we know that \text{$p\mid d$} and that the reduction of $f$ in $k$ is in the form $\fbar(z)=\Fbar\!\left(z^{p^\ell}\right)$, then we can get a better understanding of how big the $d_i$ have to be and which $x_i$ will be present in the union.

\begin{proposition}\label{prop:general-degree-preimage-of-disks}

Let $f\in\Ocal[z]$ be a monic degree $d$ polynomial and let $x\in\Ocal$. Then there exists an integer $n\geqslant 1$, integers $1\leqslant d_1,\dots,d_n\leqslant d$, and some $0<\de_x<1$, such that for all $\de_x< s<1$, 
    \[f^{-1}\!\left(\Dbar(x,s)\right)=\bigcup_{i=1}^{n}\Dbar\!\left(x_i,s^{1/d_i}\right),\]
where each $x_i\in f^{-1}(x)$. Moreover, we have that $d_1+\cdots+d_n=d$, and $f$ maps $\Dbar\!\left(x_i,s^{1/d_i}\right)$ $d_i$-to-$1$ (counting multiplicity) onto $\Dbar(x,s)$. And, if \text{$p\mid d$} and the reduction of $f$ in the form $\fbar(z)=\Fbar\!\left(z^{p^\ell}\right)$ for some nonconstant $\Fbar\in k[z]$ and $\ell\geqslant1$ (for example, any such $f$ that is postcritically bounded and fixes $0$ will work), then $p^\ell|d_i$ for all $i=1,\dots,n$. 
\end{proposition}
\begin{proof}
Let $\ze=\ze(0,1)\in\Pbb^1_\an$ be the Gauss point, let $\vvec=\vvec_\ze(x)$, and let $f_\sharp^{-1}(\vvec)=\{\wvec_1,\dots,\wvec_n\}$, where we know that each $\wvec_i$ is a direction at $\ze$ because $f^{-1}(\ze)=\ze=f(\ze)$ (since $f^{-1}(\Ocal)=\Ocal=f(\Ocal)$). Since $f$ is monic with integral coefficients, we know that $f$ has good reduction (or in the terminology of \cite{RobBook}, \emph{explicit good reduction}), and thus \cite[Theorem 7.34]{RobBook} shows that 
    \[f(\wvec_i)=f_\sharp(\wvec_i)=\vvec,\]
so $f^{-1}(x)\cap\wvec_i$ is non empty. Hence we can write $\wvec_i=\vvec_\ze(x_i)$ for some $x_i\in f^{-1}(x)\subseteq\Ocal$. Let $d_i=\deg_{\ze,\wvec_i}(f)$, and let $0<\la_{x_i}<1$ be such that $f(z)-f(x_i)=f(z)-x$ has no roots on the annulus $\Acal(x_i ; \la_{x_i},1)$. Then by definition, we have that $d_i$ is the common inner and outer Weierstrass degree of $f(z)-x$ on $\Acal(x_i; \la_{x_i},1)$. Hence Lemma \ref{lem:image-of-disc-under-monic-poly} shows that if $\la_{x_i}<s<1$, then 
    \[f\!\left(\Dbar(x_i,s)\right)=\Dbar\!\left(x,s^{d_i}\right).\]
Now let 
    \[\de_x=\max_{i=1,\dots,n}\la_{x_i}^{d_i}.\]
Then for any $\de_x< s<1$, we have that $\la_{x_i}< s^{1/d_i}<1$, and thus 
    \[f\!\left(\Dbar\!\left(x_i,s^{1/d_i}\right)\right)=\Dbar(x,s).\]
It follows that 
    \[\bigcup_{i=1}^{n}\Dbar\!\left(x_i,s^{1/d_i}\right)\subseteq f^{-1}\!\left(\Dbar(x,s)\right).\]
Since $f\!\left(\Dbar\!\left(x_i,s^{1/d_i}\right)\right)=\Dbar(x,s)$, we know that each of the $\Dbar\!\left(x_i,s^{1/d_i}\right)$ will be one of the disk components of $f^{-1}\!\left(\Dbar(x,s)\right)$. Now, if $y\in f^{-1}(x)$ then $y\in\wvec_i$ for some $i$, and thus $\vvec_\ze(y)=\vvec_\ze(x_i)$ for some $i$. But this is just to say that $|y-x_i|_v<1$, so the fact that fact that $f(z)-x$ has no roots on the annulus $\Acal(x_i ; \la_{x_i},1)$ shows that $|y-x_i|_v\leqslant\la_{x_i}<s^{1/d_i}$, and thus $y\in\Dbar\!\left(x_i,s^{1/d_i}\right)$. Hence
    \[\bigcup_{i=1}^{n}\Dbar\!\left(x_i,s^{1/d_i}\right)= f^{-1}\!\left(\Dbar(x,s)\right).\]
Furthermore, \cite[Theorems 7.30, 7.34]{RobBook} show that 
    \[\sum_{i=1}^{n}d_i=\sum_{i=1}^{n}\deg_{\ze,\wvec_i}(f)=\deg\!\left(\fbar\right)=d,\]
and \cite[Lemma 7.35]{RobBook} shows that $\deg_{\ze,\wvec_i}(f)$ is the algebraic multiplicity of $\xbar_i$ as a root of $\fbar(z)-\xbar$. Now assume that $f$ maps $\Dbar\!\left(x_i,s^{1/d_i}\right)$ $d'_i$-to-$1$ (counting multiplicity) onto $\Dbar(x,s)$. Then  $d_i$ is the algebraic multiplicity of $\xbar_i$ as a root of $\fbar(z)-\xbar$, which is the sum of the multiplicities of the roots $y\in\Ocal$ of $f(z)-x$ such that $\ybar=\xbar_i$ (since $f$ is monic with integral coefficients, and thus $f^{-1}(\Ocal)=\Ocal$), which is the sum of the multiplicities of the roots $y\in\Ocal$ of $f(z)-x$ such that $|y-x_i|_v<1$. Now, $|y-x_i|_v<1$ if and only if $|y-x_i|_v\leqslant\la_{x_i}$, since $f(z)-x$ has no roots on the anulus $\Acal(x_i; \la_{x_i},1)$. Hence $d_i$ is the sum of the multiplicities of the roots $y\in\Ocal$ of $f(z)-x$ such that $y\in\Dbar(x_i,\la_{x_i})$, and since $\la_{x_i}< s^{1/d_i}$, this shows that $d_i\leqslant d'_i$. Then Lemma \ref{lem:Robbook-theorem-3.23} shows that 
    \[d=\sum_{i=1}^{n}d_i\leqslant \sum_{i=1}^{n}d'_i=d,\]
so we must have equality throughout, meaning $d_i=d'_i$.

Now suppose that \text{$p\mid d$} and that the reduction of $f$ in $k$ is in the form $\fbar(z)=\Fbar\!\left(z^{p^\ell}\right)$ for some non-constant $\Fbar\in k[z]$ and $\ell\geqslant1$. The fact that $\Cbb_v$ is algebraically closed shows that $k$ is algebraically closed, so the assumption on the reduction type of $f$ shows that $\fbar(z)=\Fbar\!\left(z^{p^\ell}\right)=\Gbar(z)^{p^\ell}$ for some $\Fbar,\Gbar\in k[z]$. Hence the multiplicity of $\xbar_i$ as a root of $\fbar(z)-\xbar$ is divisible by $p^\ell$, so $p^\ell|\deg_{\ze,\wvec_i}(f)=d_i$.
\end{proof}

While Lemma \ref{lem:image-of-disc-under-monic-poly} and Proposition \ref{prop:general-degree-preimage-of-disks} show that the images of large enough disks under monic polynomials in $\Ocal[z]$ will behave nicely, it does nothing to establish what ``large enough" means in relation to a center of the disk. That is, Lemma \ref{lem:image-of-disc-under-monic-poly} shows that for any $x\in\Ocal$, there is some $0<\la_x<1$ and $d_x=1,\dots, d$ for which 
    \[f\!\left(\Dbar(x,s)\right)=\Dbar\!\left(f(x),s^{d_x}\right)\]
for all $\la_x<s<1$. However, it says nothing about what the smallest such value of $\la_x$ can be in relation to $x$. In fact, it is sometimes possible for the smallest such value of $\la_x$ to get arbitrarily close to $1$ as $x$ ranges over all of $\Ocal$. While we cannot get a uniform such value of $\la$ for all $x\in\Ocal$, we can get a uniform value of $\la$ if we require $x$ to not be in certain directions at the Gauss point.

\begin{proposition}\label{prop:Uniform-lambda-excluding-bad-directions}
Let $f\in\Ocal[z]$ be a monic polynomial of degree $d$. Then there is some $0<\la_0<1$ and integer $1\leqslant d_0\leqslant d$ such that for all but finitely many directions $\vvec$ (not containing $\infty$) at the Gauss point $\ze=\ze(0,1)\in\Pbb^1_\an$, if $x\in\vvec$ is a type \RomanNumeralCaps{1} point and $\la_0<s<1$, then 
    \[f\!\left(\Dbar(x,s)\right)=\Dbar\!\left(f(x),s^{d_0}\right).\]
Moreover, if \text{$p\mid d$} and the reduction of $f$ is in the form $\fbar(z)=\Fbar\!\left(z^{p^\ell}\right)$ for some nonconstant $\Fbar\in k[z]$ and $\ell\geqslant1$, then $p^\ell|d_0$. 
\end{proposition}
\begin{proof}
Write 
    \[f(z)=\sum_{i=0}^{d}a_iz^i,\]
where $a_d=1$ and $a_i\in\Ocal$. Then if $\pi\in D(0,1)$ we have that
    \[f(z+\pi)=\sum_{i=0}^{d}a_i(z+\pi)^i=a_0+\sum_{i=1}^{d}\sum_{j=0}^{i}a_i\binom{i}{j}z^{i-j}\pi^j=f(z)+\sum_{i=1}^{d}\sum_{j=1}^{i}a_i\binom{i}{j}z^{i-j}\pi^j,\]
so
    \[f(z+\pi)-f(z)=\sum_{i=1}^{d}\sum_{j=1}^{i}a_i\binom{i}{j}z^{i-j}\pi^j=\sum_{j=1}^{d}\sum_{i=j}^{d}a_i\binom{i}{j}z^{i-j}\pi^j\coloneqq\sum_{j=1}^{d}f_j(z)\pi^j,\]
where $f_j\in\Ocal[z]$ is defined so that equality holds. Note that $\deg f_j=d-j$ and that $f_d(z)=a_d=1$. Let $1\leqslant d_0\leqslant d$ be the smallest integer for which $|f_{d_0}(x)|_v=1$ for some $x\in\Ocal$. We claim that there is some $0<\la_0<1$ such that 
    \[|f(x+\pi)-f(x)|_v=|\pi|_v^{d_0}\]
for all $x\in\Ocal$ such that $|f_{d_0}(x)|_v=1$ and $\pi\in D(0,1)$ such that $|\pi|_v>\la_0$. Suppose that $x\in\Ocal$ is such that $|f_{d_0}(x)|_v=1$, and for all $1\leqslant j<d_0$, let
    \[M_j=\max_{z\in\Ocal}|f_j(z)|_v.\]
The maximum modulus principle shows that this maximum is achieved and well defined on $\Ocal$, and the fact that $f_j\in\Ocal[z]$ shows that $f_j(\Ocal)\subseteq\Ocal$, so the definition of $d_0$ shows that $M_j<1$ since $j<d_0$. Now let
    \[\la_0=\max_{1\leqslant j<d_0}M_j^{1/(d_0-j)}<1.\]
If $|\pi|_v>\la_0$ and $1\leqslant j<d_0$, then 
    \[|\pi|_v^{d_0-j}>\la_0^{d_0-j}\geqslant M_j\geqslant |f_j(x)|_v,\]
and thus 
    \[|\pi|_v^{d_0}>|f_j(x)|_v|\pi|_v^j=\left|f_j(x)\pi^j\right|_v.\]
If $d_0<j\leqslant d$, then 
    \[\left|f_j(x)\pi^j\right|_v=|f_j(x)|_v|\pi|_v^j\leqslant|\pi|_v^j<|\pi|_v^{d_0},\]
since $|f_j(x)|_v\leqslant1$ and $\pi\in D(0,1)$. Hence for all $\pi\in D(0,1)$ such that $|\pi|_v>\la_0$, we have that 
    \[|\pi|_v^{d_0}>\left|f_j(x)\pi^j\right|_v,\]
and thus 
    \[|f(x+\pi)-f(x)|_v=\left|\sum_{j=1}^{d}f_j(x)\pi^j\right|_v=|\pi|_v^{d_0}.\]
We now claim that if $x\in\Ocal$ and $|f_{d_0}(x)|_v\neq1$, then $x$ is in one of finitely many directions (not containing $\infty$) at $\ze$. In this case, we know that $|f_{d_0}(x)|_v<1$, and thus $\fbar_{d_0}(\xbar)=0$ in $k$. But $\fbar_{d_0}$ has only finitely many zeroes in $k$, say $\xbar_1,\dots,\xbar_n\in k$ for $x_1,\dots,x_n\in\Ocal$, so it follows that $\xbar=\xbar_i$ for some $i=1,\dots,n$, and thus $|x-x_i|_v<1$. But this is just to say that $x\in\vvec_\ze(x_i)$, so if $|f_{d_0}(x)|_v\neq1$, then $x\in\vvec_\ze(x_i)$ for some $i=1,\dots,n$. It follows that if $\vvec$ is a direction at $\ze$ not containing $\infty$, and $\vvec\neq\vvec_\ze(x_i)$ for all $i=1,\dots,n$, then for any type \RomanNumeralCaps{1} point $x\in\vvec$ and $\pi\in D(0,1)$ such that $|\pi|_v>\la_0$,
    \[|f(x+\pi)-f(x)|_v=|\pi|_v^{d_0}.\]
It immediately follows that 
    \[f\!\left(\Dbar(x,s)\right)=\Dbar\!\left(f(x),s^{d_0}\right)\]
for any type I point $x\in\vvec$ and $\la_0<s<1$. 

Now suppose that \text{$p\mid d$} and that the reduction of $f$ in $k$ is in the form $\fbar(z)=\Fbar\!\left(z^{p^\ell}\right)$ for some non-constant $\Fbar\in k[z]$ and $\ell\geqslant1$. Increasing $s$ if necessary so that we may apply Proposition \ref{prop:general-degree-preimage-of-disks}, we see that $d_0=d_i$ (in the notation of Proposition \ref{prop:general-degree-preimage-of-disks}) for some $i=1,\dots,n$, and thus $p^\ell|d_0$. 
\end{proof}

In the notation of the proof of Proposition \ref{prop:Uniform-lambda-excluding-bad-directions}, we will refer to the directions $\vvec_\ze(x_1),\dots,\vvec_\ze(x_n)$ as the \emph{bad directions} of $f$. As we will see in section \ref{sect:towards-infinite-wildy-ramification}, the bad directions play a key role in establishing the wild ramification of the field extensions we are interested in. 

\begin{remark}
Proposition \ref{prop:Uniform-lambda-excluding-bad-directions} becomes particularly nice if we assume that $f\in\Ocal[z]$ is a monic postcritically bounded polynomial of degree $d=p^\ell$. Then Lemma \ref{lem:wildplacegoodreduction'-composite-degree} shows that the reduction of $f$ in $k$ is in the form $\fbar(z)=z^d$. Furthermore,  
    \[\left|\binom{d}{j}\right|_v\leqslant1\]
for all $j=1,\dots,d$, with equality if and only if $j=d$, and we know that $|a_i|_v<1$ for all $i=0,\dots,d-1$. Hence $|f_j(x)|_v<1$ for all $x\in \Ocal$ if $j<d$, and $f_d(x)=1$. This shows that we have $d_0=d$ and $f$ has no bad directions, so 
    \[f\!\left(\Dbar(x,s)\right)=\Dbar\!\left(f(x),s^d\right)\]
for all $x\in \Ocal$ and $\la_0<s<1$.

\end{remark}

\section{Establishing a Lower Bound for Ramification Indices}\label{sect:Establishing-a-Lower-Bound-for-Ramification-Indices}

\hspace{\parindent} This nice behaviour of preimages of large enough discs under $f$ allows us to prove something that will be crucial in establishing a lower bound on the ramification indices of certain field extensions.

\begin{proposition}\label{prop:weaker-version-of-the-good-property}
Let $f\in\Ocal[z]$ be a monic degree $d$ (with \text{$p\mid d$}) polynomial with reduction in the form $\fbar(z)=\Fbar\left(z^{p^\ell}\right)$ for some nonconstant $\Fbar\in k[z]$ and $\ell\geqslant1$ (for example, we could have $f$ is as in Lemma \ref{lem:wildplacegoodreduction'-composite-degree}). Suppose that $\al,\be\in\Ocal$ and that $\#f^{-1}(\be)=d$. If $v(\al-\be)>0$, then for every $\al'\in f^{-1}(\al)$, there are at least $p^\ell$ values of $\be'\in f^{-1}(\be)$ for which $v(\al'-\be')>0$. 
\end{proposition}
\begin{proof}
Let $\ze=\ze(0,1)\in\Pbb^1_\an$ be the Gauss point. Then the condition that $v(\al-\be)>0$ is equivalent to the condition that $\vvec_\ze(\al)=\vvec_\ze(\be)$. If $\al'\in f^{-1}(\al)$ and $\be'\in f^{-1}(\be)$, then the fact that $\fbar$ is nonconstant and \cite[Theorem 7.34]{RobBook} shows that
\begin{equation}
    f_\sharp(\vvec_\ze(\al'))=\vvec_\ze(f(\al'))=\vvec_\ze(\al)=\vvec_\ze(\be)=\vvec_\ze(f(\be'))=f_\sharp(\vvec_\ze(\be')).\label{eq:fsharp-of-direction}
\end{equation}
Write $f_\sharp^{-1}(\vvec_\ze(\be))=\{\vvec_\ze(\be_1),\dots,\vvec_\ze(\be_n)\}$, where $\be_1,\dots,\be_n\in f^{-1}(\be)$. Then equation (\ref{eq:fsharp-of-direction}) shows that $\vvec_\ze(\al')\in f_\sharp^{-1}(\vvec_\ze(\be))$, and thus $\vvec_\ze(\al')=\vvec_\ze(\be_i)$ for some $i$. By the proof of Proposition \ref{prop:general-degree-preimage-of-disks}, we know that for $s<1$ large enough and $d_i=\deg_{\ze,\vvec_\ze(\be_i)}(f)$, we have that $f$ maps the disk $\Dbar(\be_i,s^{1/d_i})$ $d_i$-to-$1$ onto $\Dbar(\be,s)$. Hence there are exactly $d_i$ preimages (not counting multiplicity, since $\#f^{-1}(\be)=d$) $\be'$ of $\be$ in $\Dbar\!\left(\be_i,s^{1/d_i}\right)$. For each of these values of $\be'$, we know that $|\be'-\be_i|_v\leqslant s^{1/d_i}<1$, and thus $\vvec_\ze(\al')=\vvec_\ze(\be_i)=\vvec_\ze(\be')$. Hence there are at least $d_i$ values of $\be'\in f^{-1}(\be)$ for which $|\al'-\be'|_v<1$, which is equivalent to saying that $v(\al'-\be')>0$. Proposition \ref{prop:general-degree-preimage-of-disks} also shows that $p^\ell\mid d_i$, so there are at least $p^\ell$ values of $\be'\in f^{-1}(\be)$ for which $v(\al'-\be')>0$. 
\end{proof}

Using this, we can establish a strict lower bound on the ramification indices of the field extensions we are interested in.

\begin{proposition}\label{prop:e_n-bound-in-composite-case}
Assume that $f\in K[z]$ is a degree $d$ monic postcritically bounded polynomial with integral coefficients that fixes $0$, and that $\ell=v_p(d)\geqslant1$. Suppose that $a\in\Kbar$ is not in the postcritical set of $f$ and that $v(a)\geqslant0$. Let $K_n=K(f^{-n}(a))$, let $K_\infty$ be the union of all the $K_n$ for $n\geqslant1$, and let $e_n=e(K_n/K)$ be the ramification index of $K_n/K$. Then there is a real constant $A>0$ and a sequence of nonnegative real numbers $(B_n)_{n\geqslant0}$, depending only on $f$ and $a$, for which $B_n/p^{\ell n}\to0$ as $n\to\infty$ and
\begin{equation*}
    e_n>\frac{p^{\ell n}}{An+B_n}.
\end{equation*}
In particular, $K_\infty/K$ is infinitely ramified. 
\end{proposition}
\begin{proof}
Note that the assumption that $a$ is not in the postcritical set of $f$ shows that $f^n(z)-a$ has $d^n$ distinct roots for all $n\geqslant0$. Thus Proposition \ref{prop:weaker-version-of-the-good-property} shows that for any $\al,\be\in T_\infty(f,a)$, if $v(\al-\be)>0$ then for all $\al'\in f^{-1}(\al)$ there are at least $p^\ell$ values of $\be'\in f^{-1}(\be)$ for which $v(\al'-\be')>0$. For $n\geqslant1$, let $M_n$ denote the number of pairs $(\al,\be)\in f^{-n}(a)\times f^{-n}(a)$ for which $\al\neq\be$ and $v(\al-\be)>0$. Then if $n\geqslant1$, we have that
\begin{align*}
    M_{n+1}&=\#\left\{\!\left.(\al',\be')\in f^{-(n+1)}(a)\times f^{-(n+1)}(a) \,\right| \al'\neq\be' \text{ and }v(\al'-\be')>0\right\}\\
    &\geqslant\sum_{\substack{(\al,\be)\in f^{-n}(a)\times f^{-n}(a) \\ \al\neq\be}}\#\left\{\!\left.(\al',\be')\in f^{-1}(\al)\times f^{-1}(\be) \,\right| v(\al'-\be')>0\right\}.
\end{align*}
If $f(\al')=\al\neq \be=f(\be')$, then we know that $v(\al'-\be')>0$ implies that 
    \[v(\al-\be)=v(f(\al')-f(\be'))\geqslant v(\al'-\be')>0,\]
so
\begin{align*}
    M_{n+1} &\geqslant \sum_{\substack{(\al,\be) \in f^{-n}(a)\times f^{-n}(a) \\ \al\neq\be, \  v(\al-\be)>0}}\#\left\{\left.(\al',\be')\in f^{-1}(\al)\times f^{-1}(\be)\, \right| v(\al'-\be')>0\right\}\\
    &=\sum_{\substack{(\al,\be)\in f^{-n}(a)\times f^{-n}(a) \\ \al\neq\be, \  v(\al-\be)>0}}\sum_{\al'\in f^{-1}(\al)}\#\left\{(\al',\be') \left|\, \be'\in f^{-1}(\be) \text{ and } v(\al'-\be')>0\right.\right\}\\
    &\geqslant\sum_{\substack{(\al,\be)\in f^{-n}(a)\times f^{-n}(a) \\ \al\neq\be, \  v(\al-\be)>0}}\sum_{\al'\in f^{-1}(\al)}p^\ell\\
    &=\sum_{\substack{(\al,\be)\in f^{-n}(a)\times f^{-n}(a) \\ \al\neq\be, \  v(\al-\be)>0}}dp^\ell\\
    &=dp^\ell M_n.
\end{align*}
Thus 
    \[M_n\geqslant dp^\ell M_{n-1}\geqslant d^2p^{2\ell}M_{n-2}\geqslant\cdots\geqslant d^{n-1}p^{\ell(n-1)}M_1\]
for all $n\geqslant1$. Note that Lemma \ref{lem:valuation-of-discriminant} and the remarks after Lemma \ref{lem:wildplacegoodreduction'-composite-degree} show that $v(\Disc(f(z)-a))>0$, and thus $M_1\neq0$. It follows that if $n\geqslant 1$, then
    \[\frac{d^{n-1}p^{\ell(n-1)}M_1}{e_n}\leqslant \frac{M_n}{e_n}\leqslant \sum_{\substack{\al,\be\in f^{-n}(a) \\ \al\neq\be}}v(\al-\be)=v(\Disc(f^n(z)-a)).\]
Here, the second inequality follows from the fact that $f$ is monic and has coefficients in $\Ocal_K$ by assumption, and thus each $\al-\be\in \Ocal_{K_n}$ for $\al,\be\in f^{-n}(a)$ since $v(a)\geqslant0$, so if $v(\al-\be)\neq0$ then $v(\al-\be)=v_{\pfrak_n}(\al-\be)/e_n\geqslant1/e_n$ for $\pfrak_n$ the prime of $K_n$. Applying Lemma \ref{lem:valuation-of-discriminant} to $g(z)=z-a$, and using the fact that $v(\Disc(z-a))=v(1)=0$, then shows that
\begin{align*}
    \frac{d^{n-1}p^{\ell(n-1)}M_1}{e_n}\leqslant v(\Disc(f^n(z)-a))&=nd^nv(d)+\sum_{i=0}^{n-1}d^iN_{z-a,n-i}\\
    &<nd^nv(d)+d^n\sum_{i=0}^{n-1}N_{z-a,n-i}=nd^nv(d)+d^nC_n,
\end{align*}
where 
    \[C_n=\sum_{i=1}^{n}N_{z-a,i}=\sum_{i=1}^{n}\sum_{j=1}^{d-1}v(f^i(\ga_j)-a).\]
Note that $C_n/p^{\ell n}\to0$ as $n\to\infty$ by Lemmas \ref{lem:-limit-of-v(fn(ga)-a)} and \ref{lem:Analysis-lemma} since $a$ is not in the postcritical set of $f$ by assumption, and $C_n\geqslant0$ since $v(\ga_j),v(a)\geqslant0$. 
Hence
\begin{equation*}
    e_n>\frac{d^{n-1}p^{\ell(n-1)}M_1}{d^n(nv(d)+C_n)}=\frac{p^{\ell n}}{An+B_n},
\end{equation*}
where $A=dp^\ell v(d)/M_1>0$ and $B_n=dp^\ell C_n/M_1\geqslant0$.
\end{proof}

Proposition \ref{prop:e_n-bound-in-composite-case} allows us to see that as $n\to\infty$, we can make the distance between certain elements of $f^{-n}(a)$ get arbitrarily close to $1$.

\begin{lemma}\label{lem:distance-between-elements-in-tree-arbitrarily-close-to-1}
Assume that $f\in K[z]$ is a degree $d$ (with \text{$p\mid d$}) monic postcritically bounded polynomial with integral coefficients that fixes $0$. Suppose that $a\in\Kbar$ is such that $v(a)\geqslant0$ and $a$ is not in the postcritical set of $f$. Let $K_n=K(f^{-n}(a))$, let $K_\infty$ be the union of all the $K_n$ for $n\geqslant1$, and let $e_n=e(K_n/K)$ be the ramification index of $K_n/K$. Then for any $0<\la<1$ and $n\geqslant1$ big enough, there are some distinct $\al_n,\be_n\in f^{-n}(a)$ such that $1>|\al_n-\be_n|_v>\la$.
\end{lemma}
\begin{proof}
With notation as in the proof of Proposition \ref{prop:e_n-bound-in-composite-case}, we have that if 
    \[V_n=\min_{\substack{\al,\be\in f^{-n}(a) \\ 0<v(\al-\be)<\infty}}v(\al'-\be')>0,\]
then 
    \[d^{n-1}p^{\ell(n-1)}M_1V_n\leqslant M_nV_n\leqslant \sum_{\substack{\al,\be\in f^{-n}(a) \\ \al\neq\be}}v(\al-\be)<d^n(nv(d)+C_n).\]
Hence
    \[V_n\leqslant \frac{d(nv(d)+C_n)}{p^{\ell(n-1)}M_1},\]
so Lemma \ref{lem:-limit-of-v(fn(ga)-a)} shows that $V_n\to0$ as $n\to\infty$ because $\ell>0$. By definition, for all $n\geqslant 1$ there are some distinct $\al_n,\be_n\in f^{-n}(a)$ such that $v(\al_n-\be_n)\neq0$ and $v(\al_n-\be_n)=V_n$. Hence $0<v(\al_n-\be_n)=V_n\to0$ as $n\to\infty$, so $1>|\al_n-\be_n|_v\to1$ as $n\to\infty$. Thus for any $\la<1$, we have that $1>|\al_n-\be_n|_v>\la$ for all sufficiently large $n$. 
\end{proof}

We are now ready to begin to relate the results we proved about images of disks under polynomials to the ramification of the field extensions we are interested in.

\begin{proposition}\label{prop:no-bad-directions-implies-wildly-ramified}
Assume that $f\in K[z]$ is a degree $d$ (with \text{$p\mid d$}) monic postcritically bounded polynomial with integral coefficients that fixes $0$. Suppose that $a\in\Kbar$ is not in the postcritical set of $f$ and that $v(a)\geqslant0$. Let $K_n=K(f^{-n}(a))$, let $K_\infty$ be the union of all the $K_n$ for $n\geqslant1$, and let $e_n=e(K_n/K)$ be the ramification index of $K_n/K$. Furthermore, suppose that there is some $0<\la_0<1$ and integer $1\leqslant d_0\leqslant d$, with \text{$p\mid d_0$}, such that for all $\la_0<s<1$ and $\al\in T_\infty(f,a)$, we have that
    \[f\!\left(\Dbar(\al,s)\right)=\Dbar\!\left(f(\al),s^{d_0}\right).\]
Then $K_\infty/K$ is infinitely wildly ramified.
\end{proposition}
\begin{proof}
By Lemma \ref{lem:distance-between-elements-in-tree-arbitrarily-close-to-1}, we can find some $N\geqslant1$ and distinct $\al,\be\in f^{-N}(a)$ for which $1>|\al-\be|_v>\la_0$. Then 
    \[1>|\al-\be|_v^{1/d_0}>|\al-\be|_v>\la_0,\]
so if $\al_1\in f^{-1}(\al)$ then 
    \[f\!\left(\Dbar\!\left(\al_1,|\al-\be|_v^{1/d_0}\right)\right)=\Dbar(\al,|\al-\be|_v).\]
Thus if $\be_1\in f^{-1}(\be)\cap\Dbar\!\left(\al_1,|\al-\be|_v^{1/d_0}\right)$ (such a $\be_1$ exists since $\be\in\Dbar(\al,|\al-\be|_v)$), then 
    \[|\al_1-\be_1|_v=|\al-\be|_v^{1/d_0}>\la_0\]
by Lemma \ref{lem:preimage-of boundary-is-on-boundary}, since $|f(\al_1)-f(\be_1)|_v=|\al-\be|_v$. Continuing in this manner, we can find sequences $(\al_n)_{n\geqslant0}$ and $(\be_n)_{n\geqslant0}$ of elements of $K_{N+n}$ for which $\al_0=\al$, $\be_0=\be$, and $f(\al_n)=\al_{n-1}$, $f(\be_n)=\be_{n-1}$, and
    \[1>|\al_n-\be_n|_v=|\al_{n-1}-\be_{n-1}|_v^{1/d_0}=\cdots=|\al-\be|_v^{1/d_0^n}>\la_0.\]
Hence 
    \[v(\al_n-\be_n)=\frac{v(\al-\be)}{d_0^n}>0.\]
If $\pfrak_n$ denotes the prime of $K_n$, so that $v|_{K_n}=v_{\pfrak_n}/e_n$, then this equation shows that 
    \[\frac{v_{\pfrak_{N+n}}(\al_n-\be_n)}{e_{N+n}}=\frac{v(\al-\be)}{d_0^n},\]
and thus
    \[e_{N+n}=\frac{v_{\pfrak_{N+n}}(\al_n-\be_n)}{v(\al-\be)}d_0^n.\]
The fact that \text{$p\mid d_0$} then shows that $p|e_{N+n}$ for all large enough $n$. Thus $K_{n}/K$ is wildly ramified for all large enough $n$, and $v_p(e_n)\to\infty$ as $n\to\infty$. 
\end{proof}

\begin{remark}
Much of what we've shown so far becomes greatly simplified in the case where $d=p^\ell$ is the power of a prime. Since $f$ will have no bad directions in this case, Proposition \ref{prop:weaker-version-of-the-good-property} shows that if $\al,\be\in\Ocal$, if $\#f^{-1}(\be)=d$, and if $v(\al-\be)>0$, then $v(\al'-\be')>0$ for all $\al'\in f^{-1}(\al)$ and $\be'\in f^{-1}(\be)$. In particular, this shows that if $a\in \Kbar$ is such that $v(a)\geqslant0$ and $a$ is not in the postcritical set of $f$, then $v(\al-\be)>0$ for all $\al,\be\in f^{-n}(a)$. 
\end{remark}

\section{Towards Infinite Wild Ramification}\label{sect:towards-infinite-wildy-ramification}

\hspace{\parindent} Suppose the normal setup, so $f\in K[z]$ is a monic degree $d$ (with \text{$p\mid d$}) postcritically bounded polynomial with integral coefficients that fixes $0$. In particular, Lemma \ref{lem:wildplacegoodreduction'-composite-degree} shows that $f$ has reduction in the form $\fbar(z)=\Fbar\!\left(z^{p^\ell}\right)$, where $\ell=v_p(d)$. Furthermore, we assume that $a\in K$ is not in the postcritical set of $f$ and that $v(a)\geqslant0$. In this section, we show that the field extensions $K(f^{-n}(a))/K$ are infinitely wildly ramified. If $d=p^\ell$ then $f$ has no bad directions, and thus Proposition \ref{prop:no-bad-directions-implies-wildly-ramified} shows that $K_\infty/K$ is infinitely wildly ramified. We thus assume that $d$ is not a power of a prime. 

Given $a\in K$ such that $v(a)\geqslant0$, we may form the iterated preimage set of $a$ under $f$, $T_\infty(f,a)$. Reducing modulo $\mfrak=D(0,1)$, we get the reduced iterated preimage set of $\abar$ under $\fbar$, $T_\infty(\fbar,\abar)$, which is also the image of $T_\infty(f,a)$ under the reduction map $\Ocal\to k$. We wish to apply the lemmas of appendix \ref{sect:some-lemmas-about-preimage-sets} to $T_\infty(\fbar,\abar)$, for which we must verify that there is at most one value of $\albar\in k$ for which $\#\fbar^{-1}(\albar)=1$. Consider the reduction $\fbar\in k[z]$ as a map of smooth curves $\fbar:\Pbb^1_k\to\Pbb^1_k$. Then when we factor $\fbar$ as $\fbar=\Fbar\circ\vp$, where $\vp(z)=z^{p^\ell}$, we know that the map $\Fbar:\Pbb^1_k\to\Pbb^1_k$ is separable, for if it was not then \cite[corollary 2.12]{ElipCurves} would show that we can factor it further as a composition of maps $\Fbar(z)=\Gbar\!\left(z^{p^r}\right)$, where $\Gbar:\Pbb^1_k\to\Pbb^1_k$ is separable and $r\geqslant1$. But $p^\ell$ is the highest power of $p$ dividing $d=\deg\fbar$, so $p\nmid\deg\Fbar$, and thus this is not possible. The Riemann--Hurwitz formula then shows that 
    \[\sum_{P\in\Pbb^1_k(k)}(e_{\Fbar}(P)-1)\leqslant2(\deg\Fbar-1). \]
Because $\infty$ is a totally ramified point of $\Fbar$ (since $\Fbar$ is a polynomial), meaning $e_{\Fbar}(\infty)=\deg\Fbar$, this implies that 
    \[\sum_{\albar\in k}(e_{\Fbar}(\albar)-1)\leqslant\deg\Fbar-1. \]
Hence there can be at most one element $\albar\in k$ for which $e_{\Fbar}(\albar)=\deg\Fbar$ (where we are using the fact that $d$ is not a power of a prime, so $\deg\Fbar>1$), which is equivalent to saying that there can be at most one point $\albar\in k$ for which $\#\Fbar{}^{-1}(\albar)=1$. But $\fbar^{-1}(\albar)=\left\{\bebar\right\}$ for some $\bebar\in k$ if and only $\Fbar{}^{-1}(\albar)=\left\{\vp\!\left(\bebar\right)\right\}$ since $\vp$ is a bijection, so there can be at most one element $\albar\in k$ for which $\#\fbar^{-1}(\albar)=1$. Applying Lemma \ref{lem:preim-tree-either-has-size-1-or-is-infinite} then shows that either $T_\infty(\fbar,\abar)$ has a single element, or $T_\infty(\fbar,\abar)$ has infinitely many elements. Now, elements of $k$ uniquely represent directions at the Gauss point $\ze=\ze(0,1)\in\Pbb^1_\mathrm{an}$ away from infinity, for if $x,y\in\Ocal$ then $\vvec_\ze(x)=\vvec_\ze(y)$ if and only if $|x-y|_v<1$, which is true if and only if $\xbar=\ybar$ in $k$. We are now ready to prove infinite wild ramification in the easier cases, which occur when we can find elements of $T_\infty(\fbar,\abar)$ that do not represent bad directions of $f$.

\begin{proposition}\label{prop:no-bad-directions-in-preimage-tree}
Assume that $f\in K[z]$ is a monic degree $d$ (with \text{$p\mid d$}) postcritically bounded polynomial with integral coefficients that fixes $0$. Suppose that $a\in K$ is not in the postcritical set of $f$ and $v(a)\geqslant0$. Let $K_n=K(f^{-n}(a))$ and let $K_\infty$ be the union of all the $K_n$ for $n\geqslant0$. Furthermore, suppose that the reduced iterated preimage set $T_\infty(\fbar,\abar)$ either has one element, and that element does not represent a bad direction of $f$, or it has infinitely many elements. Then $K_\infty/K$ is infinitely wildy ramified. 
\end{proposition}
\begin{proof}
First suppose that $T_\infty(\fbar,\abar)$ has one element, and that element does not represent a bad direction of $f$. Then no element of $T_\infty(f,a)$ will lie in a bad direction of $f$, so it follows by the definition of the bad directions of $f$ that there is some $0<\la_0<1$ and integer $1\leqslant d_0\leqslant d$, with \text{$p\mid d_0$}, such for all $\al\in T_\infty(f,a)$ and $\la_0<s<1$, we have that
    \[f\!\left(\Dbar(\al,s)\right)=\Dbar\!\left(f(\al),s^{d_0}\right).\]
Proposition \ref{prop:no-bad-directions-implies-wildly-ramified} then shows that $K_\infty/K$ is infinitely wildly ramified. 

Now suppose that $T_\infty(\fbar,\abar)$ has infinitely many elements. Then Lemma \ref{lem:avoiding-elements-in-infinite-tree}\ref{avoiding-elements-in-infinite-tree-b} shows that we can find some $a_0\in T_\infty(f,a)$ for which $T_\infty(\fbar,\abar_0)$ contains no elements representing one of the finitely many bad directions of $f$. Hence $T_\infty(f,a_0)$ contains no elements that lie in a bad direction of $f$. We still have that $v(a_0)\geqslant0$ because $f^{-1}(\Ocal)\subseteq\Ocal$, and $a_0$ is not in the postcritical set of $f$ since $a_0\in T_{\infty}(f,a)$ and $a$ is not in the postcritical set of $f$, so $f^n(z)-a_0$ has $d^n$ distinct roots for all $n\geqslant0$. If $L_n=K(f^{-n}(a_0))$, then Proposition \ref{prop:no-bad-directions-implies-wildly-ramified} shows that $L_n/K$ is wildy ramified (for large enough $n$) with $v_p(e(L_n/K))$ diverging to $\infty$, so the same is true for $K_n/K$. Thus $K_\infty/K$ is infinitely wildly ramified. 
\end{proof}

\subsection{A Totally Invariant Bad Direction}

\hspace{\parindent} It follows by Lemma \ref{lem:preim-tree-either-has-size-1-or-is-infinite} that the only situation not dealt with in Proposition \ref{prop:no-bad-directions-in-preimage-tree} is the one in which $T_\infty(\fbar,\abar)$ contains one element, and that element represents a bad direction of $f$. In this case, we have that $\vvec_\ze(a)$ is a bad direction of $f$, and that 
    \[f^{-1}(\vvec_\ze(a))=\vvec_\ze(a)=f(\vvec_\ze(a)),\] 
or put in terms of disks, that
    \[f^{-1}(D(a,1))=D(a,1)=f(D(a,1)).\]
We first claim that there is a $K$-conjugate of $f$ with monomial good reduction. 

\begin{lemma}\label{lem:conjugate-with monomial-good-reduction}
Let $f,a$ be as above. Then there exists a monic $K$-conjugate $g$ of $f$ with integral coefficients for which $\gbar(z)=z^d$.  
\end{lemma}
\begin{proof}
Define $\mu\in K[z]$ by $\mu(z)=z-a$, and let 
    \[g(z)\coloneqq \mu\circ f\circ\mu^{-1}(z)=f(z+a)-a,\]
so $\vvec_\ze(0)=\vvec_\ze(\mu(a))$ is a totally invariant direction of $g$. We claim that $g$ has monomial good reduction, so if 
    \[g(z)=z^d+\sum_{i=0}^{d-1}c_iz^i,\]
then we claim that $v(c_i)>0$ for all $i=0,\dots,d-1$. The fact that $v(a)\geqslant0$ shows that $v(c_i)\geqslant0$ for all $i=0,\dots,d-1$. Note that 
    \[v(c_0)=v(g(0))=v(f(a)-a)>0,\]
so suppose that $v(c_i)=0$ for some $i=1,\dots,d-1$. Then for $c_0\neq c\in D(0,1)$, the Newton polygon of $g(z)-c$ will be the lower convex hull of the points 
    \[(0,v(c_0-c)), (1,v(c_1)),\dots (i-1,v(c_{i-1})), (i,0), (i+1, v(c_{i+1})), \dots (d-1, v(c_{d-1})), (d,0).\]
We have that $v(c_j)\geqslant0$ for all $j=0,\dots,d-1$, so this shows that the Newton polygon of $g(z)-c$ will have a line segment of slope $0$, and thus $g(z)-c$ has a root of absolute value $1$. But $c\in D(0,1)$, so the fact that $\vvec_\ze(0)$ is completely invariant under $g$ shows that $g(z)-c$ must have all of its roots in $D(0,1)$, so it cannot have a root of absolute value $1$. Hence $v(c_i)>0$ for all $i=0,\dots,d-1$, so $g$ has monomial good reduction.
\end{proof}

Because $\mu\in K[z]$ and $\mu(a)=0$, we know that 
    \[K(g^{-n}(0))=K(\mu(f^{-n}(a)))=K(f^{-n}(a)).\]
Hence the the claim that $K(f^{-n}(a))/K$ is wildly ramified is equivalent to the claim that $K(g^{-n}(0))/K$ is wildly ramified.

\begin{proposition}\label{prop:one-totally-invariant-bad-direction-nice-case}
Assume the setup as in Lemma \ref{lem:conjugate-with monomial-good-reduction}. Then for any $\be_0\in\Kbar\cap D(0,1)$ such that $v(\be_0)\neq v(c_0)$, the extensions $F_n=K(g^{-n}(\be_0))$ are wildly ramified over $K$ for sufficiently large $n$, and $v_p(e(F_n/K))\to\infty$ as $n\to\infty$.

\end{proposition}
\begin{proof}
First suppose that $v(\be_0)<v(c_0)$. Then we have that $g(z)-\be_0$ has constant term $c_0-\be_0$ with valuation $v(c_0-\be_0)=v(\be_0)$, so the Newton polygon $N(g(z)-\be_0)$ of $g(z)-\be_0$ will be the lower convex hull of the points $(0,v(\be_0)),(1,v(c_1)),\dots,(d-1,v(c_{d-1})), (d,0)$. Hence there will be a root $\be_1\in F_1$ of $g(z)-\be_0$ with valuation
    \[v(\be_1)=-\frac{v(\be_0)-v(c_i)}{0-i}=\frac{v(\be_0)-v(c_i)}{i}<v(\be_0)<v(c_0)\]
for some $i=1,\dots,d$ (where we are taking $c_d=1$). Note that $v(\be_1)>0$ by the assumption that $\be_0\in D(0,1)$ and that $\vvec_\ze(0)$ is totally invariant under $g$. Continuing in this manner, we see that we can find a sequence of elements $\be_n\in F_n$ such that $g(\be_n)=\be_{n-1}$ and such that for each $n$, there is some $i_n=1,\dots,d$ for which
    \[0<v(\be_n)=\frac{v(\be_{n-1})-v(c_{i_n})}{i_n}<v(\be_{n-1})<\cdots<v(\be_0)<v(c_0).\]
Note that $i_n$ can be $1$ for only a finite number of values of $n$, since if $i_n=1$ then $v(\be_n)=v(\be_{n-1})-v(c_1)$, and $v(\be_n)>0$ for all $n$. Hence for sufficiently large $n$, we have that
    \[v(\be_n)=\frac{v(\be_{n-1})-v(c_{i_n})}{i_n}<\frac{v(\be_{n-1})}{i_n}\leqslant \frac{v(\be_{n-1})}{2}. \]
It follows that for all sufficiently large $n$, we will have that $v(\be_n)<v(c_i)$ for all $i=0,\dots,d-1$. Thus the Newton polygon $N(g(z)-\be_n)$ will contain only one segment of slope $-v(\be_n)/d$, which means that for any $\be_{n+1}\in g^{-1}(\be_n)$, 
    \[v(\be_{n+1})=\frac{v(\be_n)}{d}.\]
Therefore, if we choose $\be_n\in F_n$ appropriately, there is some $N\geqslant0$ such that if $n\geqslant 0$ then 
    \[v(\be_{N+n})=\frac{v(\be_N)}{d^n}.\]
If we write $\pfrak_n$ for the prime of $F_n$, so the extension of $v$ to $F_n$ is given by $v|_{F_n}=v_{\pfrak_n}/e(F_n/K)$, then this shows that 
    \[e(F_{N+n}/K)=\frac{v_{\pfrak_N+n}(\be_{N+n})}{v(\be_N)}d^n,\]
and thus $v_p(e(F_n/K))\to\infty$ as $n\to\infty$. 

Now suppose that $v(\be_0)>v(c_0)$ (so in particular, $v(c_0)<\infty$). Then $g(z)-\be_0$ has constant term $c_0-\be_0$ with valuation $v(c_0-\be_0)=v(c_0)$, and thus the Newton polygon $N(g(z)-\be_0)$ is given by the lower convex hull of the points 
    \[(0,v(c_0)), (1,v(c_1)), \dots, (d-1,v(c_{d-1})), (d,0).\] 
Hence $g(z)-\be_0$ has a root $\be_1\in F_1$ with valuation
    \[0<v(\be_1)=\frac{v(c_0)-v(c_i)}{i}<v(c_0)\]
for some $i=1,\dots,d$. Applying the first part of this proof to a sequence starting at $\be_1$ then shows that $v_p(e(F_n/K))\to\infty$ as $n\to\infty$. 
\end{proof}

\begin{proposition}\label{prop:totally-invariant-bad-direction-wildly-ramified}
Assume the setup as in Lemma \ref{lem:conjugate-with monomial-good-reduction}. Then $K(g^{-n}(0))/K$ is wildly ramified over $K$ for large enough $n$ and $v_p(e(K(g^{-n}(0))/K))\to\infty$ as $n\to\infty$. 
\end{proposition}
\begin{proof}
If $\infty=v(0)\neq v(c_0)$, then Proposition \ref{prop:one-totally-invariant-bad-direction-nice-case} gives the desired result. Now suppose that $v(c_0)=\infty$, so $c_0=0$. The assumption that $\#T_\infty(f,a)=\infty$, and thus $\#T_\infty(g,0)=\infty$, shows that we cannot have $v(\be_0)=\infty$ for all $\be_0\in T_\infty(g,0)$, so we must have $v(\be_0)<\infty$ for some $\be_0\in T_\infty(g,0)$. Applying Proposition \ref{prop:one-totally-invariant-bad-direction-nice-case} to $\be_0$ then shows that $v_p(e(K(g^{-n}(\be_0))/K))\to\infty$, so we also have that $v_p(e(K(g^{-n}(0))/K))\to\infty$.
\end{proof}

\InfinitelyRamifiedLocalField*
\begin{proof}
If $d=p^\ell$ is a prime power, then $f$ has no bad directions, so $K_\infty/K$ is infinitely wildly ramified by Proposition \ref{prop:no-bad-directions-implies-wildly-ramified}. Now assume that $d$ is not a power of a prime. Then we know by Lemma \ref{lem:wildplacegoodreduction'-composite-degree} that $f$ has reduction in the form $\fbar(z)=\Fbar\!\left(z^{p^\ell}\right)$, where $\ell=v_p(d)\geqslant1$ and $\Fbar\in k[z]$ is nonconstant. Then Proposition \ref{prop:no-bad-directions-implies-wildly-ramified} shows that $K_\infty/K$ is infinitely wildly ramified if one of the following holds: we have $\#T_\infty(\fbar,\abar)=1$, and its one element does not lie in a bad direction of $f$; or we have that $T_\infty(\fbar,\abar)$ is infinite. Furthermore, Proposition \ref{prop:totally-invariant-bad-direction-wildly-ramified} show that $K_\infty/K$ is infinitely wildly ramified if $T_\infty(\fbar,\abar)$ contains a single element, and that element represents a bad direction of $f$. This deals with all of the possible cases, which shows that $K_\infty/K$ is always infinitely wildly ramified. 
\end{proof}

\section{Proof of the Main Theorem}
Now that we have proved Theorem \ref{thm:InfinitelyRamifiedLocalField}, we are ready to apply it to prove Theorem \ref{thm:InfinitelyRamified}. The only things stopping us from immediately doing this are some mild additional assumptions on $f$. However, we can replace $f$ with an appropriate $\Kbar$-conjugate and replace the base point with an element in its iterated preimage set under $f$ in order to force it to satisfy the conditions necessary  to apply 
Theorem \ref{thm:InfinitelyRamifiedLocalField}.

\InfinitelyRamified*
\begin{proof}
First suppose that $a$ is not in the postcritical set of $f$. Let $K_v$ denote the completion of $K$ with respect to $v$ and $\Kbar_v$ denote a fixed algebraic closure of $K_v$. Because $f\in K_v[z]$ has potential good reduction, there is some $\mu_1(z)=\al z+\be\in\Kbar_v[z]$ for which $\mu_1\circ f\circ\mu_1^{-1}\in \Kbar_v[z]$ is monic with integral coefficients. Passing to the finite extension $K_v(\al,\be)/K_v$ if necessary, we may assume that $\mu_1\in K_v[z]$ and that $\mu_1\circ f\circ\mu_1^{-1}\in\Ocal_{K_v}[z]$. Then the roots of $\mu_1\circ f\circ\mu_1^{-1}(z)-z$ all lie in $\Ocal_{\Kbar_v}$, so we let $b\in\Ocal_{\Kbar_v}$ be a root of $\mu_1\circ f\circ\mu_1^{-1}(z)-z$, and we let $\mu_2(z)=z-b\in\Kbar_v[z]$. Passing to the extension $K_v(b)/K_v$ if necessary, we may also assume that $\mu_2\in K_v[z]$. Hence if $\mu=\mu_2\circ\mu_1$, then $g\coloneqq \mu\circ f\circ\mu^{-1}\in K_v[z]$ is a degree $d$ monic polynomial with integral coefficients that fixes $0$. The fact that $f$ is postcritically bounded clearly implies that $g$ is postcritically bounded. Let $K_{v,n}=K_v(g^{-n}(\mu(a)))=K_v(f^{-n}(a))$ and let $K_{v,\infty}$ denote the union of all the $K_{v,n}$ for $n\geqslant0$. It suffices to show that $K_{v,\infty}/K_v$ is infinitely wildly ramified in order to show that $K_\infty/K$ is infinitely wildly ramified at $\pfrak$. If $v(\mu(a))<0$, then the Newton polygon of $g(z)-\mu(a)$ will consist of a single line segment connecting $(0,v(\mu(a)))$ to $(d,0)$, with slope $-v(\mu(a))/d$, and thus every root of $g(z)-\mu(a)$ will have valuation $v(\mu(a))/d<0$. Continuing in this manner, we see that if $\al_n\in g^{-n}(\mu(a))$ then $v(\al_n)=v(\mu(a))/d^n$, which shows that $K_{v,\infty}/K_v$ is infinitely wildly ramified. We thus assume that $v(\mu(a))\geqslant0$. Then the fact that $a$ is not in the postcritical set of $f$ shows that $\mu(a)$ is not in the postcritical set of $g$. It follows that $(K_v,v)$, $g\in\Ocal_{K_v}[z]$, and $\mu(a)\in K_v   $ satisfy all the hypotheses of Theorem \ref{thm:InfinitelyRamifiedLocalField}, so the extension $K_{v,\infty}/K_v$ is infinitely wildly ramified, and consequently $K_\infty/K$ is infinitely wildly ramified at $\pfrak$.

Now suppose that $a$ is in the postcritical set of $f$. Note that if $\ga\in\Kbar$ is a critical point of $f$ and $f^n(\ga)\in T_{\infty}(f,a)$ for some $n\geqslant0$ then $\ga\in T_\infty(f,a)$. But Lemma \ref{lem:avoiding-elements-in-infinite-tree}\ref{avoiding-elements-in-infinite-tree-a} then shows that there will only be finitely many elements of the forward orbit of $\ga$ under $f$ contained in $T_\infty(f,a)$, so the fact that $f$ has finitely many critical points shows that there will only be finitely many postcritical elements of $f$ contained in $T_\infty(f,a)$. Thus Lemma \ref{lem:avoiding-elements-in-infinite-tree}\ref{avoiding-elements-in-infinite-tree-b} shows that we can find some $a_0\in T_\infty(f,a)$ for which $T_\infty(f,a_0)$ is disjoint from the postcritical set of $f$ (and thus $\#T_\infty(f,a_0)=\infty$). We can thus apply the results shown in the first part of this proof (after replacing $K$ with $K(a_0)$) to see that $K_\infty/K$ is infinitely wildly ramified above $\pfrak$. 
\end{proof}

\newpage
\appendix
\section{Appendix}
\subsection{Dealing With Limit Points of the Postcritical Set}

\hspace{\parindent} In the process of proving Proposition \ref{prop:e_n-bound-in-composite-case}, we encountered a sequence of nonnegative real numbers in the form 
    \[C_n=\sum_{i=1}^{n}\sum_{j=1}^{d-1}v(f^i(\ga_j)-a),\]
where $\ga_1,\dots,\ga_{d-1}$ are the critical points of $f$ and $a\in\Dbar(0,1)$ is not in the postcritical set of $f$. While we cannot expect $C_n$ to be bounded as $n\to\infty$, we now show that $C_n$ is dominated by any function that grows faster than $n$, which is strong enough to prove Proposition \ref{prop:e_n-bound-in-composite-case}. We first show a result about what happens when an open disk of unit radius is periodic under a polynomial satisfying certain nice conditions. The following arguments were first told to the author in a private communication with Rob Benedetto.

\begin{lemma}\label{lem:periodic-disk-has-unique-attracting-periodic-point}
Suppose that $f\in\Ocal[z]$ is a polynomial for which $f'\!\left(\Dbar(0,1)\right)\subseteq D(0,1)$, and that $V=D(a,1)\subseteq\Dbar(0,1)$ is periodic under $f$ of period $\ell\geqslant1$, so $f^\ell(V)=V$. Then $V$ contains a unique attracting periodic point for $f$. That is, there is a unique periodic point $P\in V$ of $f$ of period $\ell$ for which
    \[\lim_{n\to\infty}f^{\ell n}(x)=P\]
for all $x\in V$. 
\end{lemma}
\begin{proof}
Let $g=f^\ell\in\Ocal[z]$, so $g(V)=V$. We first claim that $V$ contains a fixed point for $g$. It will suffice to show that if $z_0\in V$ then 
    \[\lim_{n\to\infty}g^n(z_0)\]
exists, for if it does then it will be a fixed point of $g$ that is contained in $V$ (since the fact that $g(V)=V$ shows that $g^n(z_0)\in V$, and thus this converges to an element of $V$ since $V$ is (topologically) closed in $\Cbb_v$). In order to show this limit exists it suffices to show that $(g^n(z_0))_{n\geqslant0}$ is a Cauchy sequence since $\Cbb_v$ is complete, and the non-archimedean triangle inequality shows that it suffices to show that $|g^{n+1}(z_0)-g^n(z_0)|_v\to0$ as $n\to\infty$. 
If $c\in V$ then we can write
    \[g(z)-g(c)=\sum_{i=1}^{d}b_i(z-c)^i,\]
where $b_1=g'(c)\in D(0,1)$ and $b_i\in\Dbar(0,1)$. Note that if $i\geqslant2$ and $x\in V$ (so $|x-c|_v<1$) then 
    \[|b_i|_v|x-c|_v^i\leqslant |b_i|_v|x-c|_v^2\leqslant |x-c|_v^2,\]
so 
    \begin{align*}
    |g(x)-g(c)|_v&\leqslant\max\!\left\{|g'(c)|_v|x-c|_v,|b_2|_v|x-c|_v^2,\dots,|b_d|_v|x-c|_v^d\right\}\\
    &\leqslant\max\!\left\{|g'(c)|_v|x-c|_v,|x-c|_v^2\right\}.
    \end{align*}
Now, note that if $x\in\Dbar(0,1)$ then 
    \[|g'(x)|_v=\left|\left(f^\ell\right)'\!(x)\right|_v=\prod_{i=0}^{\ell-1}\left|f'(f^i(x))\right|_v<1,\]
so $g'\!\left(\Dbar(0,1)\right)\subseteq D(0,1)$, and thus $g'\!\left(\Dbar(0,1)\right)=\Dbar(\al,r)$ for some $\al\in D(0,1)$ and $r<1$. Hence there is some constant $0<R<1$, where $R=\max\{|\al|_v,r\}$, for which $g'\!\left(\Dbar(0,1)\right)\subseteq\Dbar(0,R)$, so $|g'(x)|_v\leqslant R<1$ for all $x\in\Dbar(0,1)$, so 
    \[|g(x)-g(c)|_v\leqslant\max\!\left\{R|x-c|_v,|x-c|_v^2\right\}.\]
Now, the fact that $g(V)=V$ is an open disk of radius $1$ shows that if $z_0\in V$ then 
\[
0<\de\coloneqq|g(z_0)-z_0|_v<1,
\]
where we assume $\de>0$, since if $\de=0$, then $z_0$ is a fixed point of $g$, which is what we are trying to show exists. Then 
\begin{align*}
    \left|g^2(z_0)-g(z_0)\right|_v
    &=|g(g(z_0))-g(z_0)|_v\\
    &\leqslant \max\!\left\{R|g(z_0)-z_0|_v,|g(z_0)-z_0|_v^2\right\}\\
    &=\max\!\left\{R\de,\de^2\right\},
\end{align*}
so
\begin{align*}
    \left|g^3(z_0)-g^2(z_0)\right|_v
    &=\left|g(g^2(z_0))-g(g(z_0))\right|_v\\
    &\leqslant\max\Bigl\{R|g^2(z_0)-g(z_0)|_v,|g^2(z_0)-g(z_0)|_v^2\Bigr\}\\
    &\leqslant R^2\de\quad\text{or}\quad R^2\de^2\quad\text{or}\quad R\de^2\quad\text{or}\quad \de^4.
\end{align*}
Continuing in this manner, we see that there are sequences of nonnegative integers $(s_n)_{n\geqslant0},(t_n)_{n\geqslant0}$ for which 
    \[|g^{n+1}(z_0)-g^{n}(z_0)|_v\leqslant R^{s_n}\de^{t_n},\]
and for which at least one of $s_n,t_n$ diverges to $\infty$ as $n\to\infty$. The fact that $0< R,\de<1$ then shows that $|g^{n+1}(z_0)-g^n(z_0)|_v\to0$ as $n\to\infty$, so $(g^n(z_0))_{n\geqslant0}$ is a Cauchy sequence in $\Cbb_v$. Hence 
    \[\lim_{n\to\infty}g^n(z_0)\in V\]
exists for any $z_0\in V$, so $g$ has a fixed point in $V$. Call one of these fixed points $P\in V$, so $g(P)=P$. Then there are constants $b_1,\dots,b_d\in\Dbar(0,1)$ for which $|b_1|_v=|g'(P)|_v<1$ and
    \[g(z)-P=g(z)-g(P)=\sum_{i=1}^{d}b_i(z-P)^i.\] 
If we define $G(z)=g(z+P)-P$, then we have that 
    \[G(z)=\sum_{i=1}^{d}b_iz^i\]
for $b_i\in\Dbar(0,1)$ with $|b_1|_v<1$, so \cite[Proposition 6.2.1]{BenedettoThesis} shows that for any $x\in D(0,1)$, we have
    \[\lim_{n\to\infty}G^n(x)=0.\]
But $G^n(z)=g^n(z+P)-P$, so this shows that 
    \[\lim_{n\to\infty}g^{n}(x+P)=P\]
for any $x\in D(0,1)$, and thus 
    \[\lim_{n\to\infty}g^n(x)=P\]
for any $x\in D(P,1)=D(a,1)=V$. From this, it follows that $P\in V$ is the unique fixed point of $g$ in $V$, and that $P$ attracts every element of $V$. The fact that $g=f^\ell$ then shows that $P$ is the unique periodic point of $f$ in $V$ with period $\ell$, and that 
    \[\lim_{n\to\infty}f^{\ell n}(x)=P\]
for any $x\in V$. 
\end{proof}

Using Lemma \ref{lem:periodic-disk-has-unique-attracting-periodic-point}, we can show that $v(f^n(\ga)-a)$ is dominated by any function that grows faster than $n$, which will be enough to show that the aforementioned sequence $(C_n)_{n\geqslant0}$ is also dominated by any such function.

\begin{lemma}\label{lem:-limit-of-v(fn(ga)-a)}
Suppose that $f\in K[z]$ is a monic postcritically bounded polynomial of degree $d$ with integral coefficients that fixes $0$. Assume that $\ell=v_p(d)\geqslant1$, that $a\in\Kbar$ with $v(a)\geqslant0$, and that $a$ does not lie in the postcritical set of $f$. Let $\psi:\Zbb\to(0,\infty)$ be any function that grows faster then $n$, i.e., $n/\psi(n)\to 0$ as $n\to\infty$. Then if $\ga\in\Kbar$ is a critical point of $f$, we have that
\begin{equation}
    \lim_{n\to\infty}\frac{v(f^n(\ga)-a)}{\psi(n)}=0.\label{eq:desired-limit}
\end{equation}
\end{lemma}
\begin{proof}
Note that remark \ref{critical-points-in-unit-disk} shows that if $\ga$ is a critical point of $f$ then $|\ga|_v\leqslant1$. We claim that $f'$ has the form $f'(z)=p^\ell g(z)$ for some $g\in\Ocal[z]$, for which it will suffice to show that if
    \[f(z)=z^d+\sum_{i=1}^{d-1}a_iz^i,\]
so 
    \[f'(z)=dz^{d-1}+\sum_{i=1}^{d-1}ia_iz^{i-1},\]
then \text{$p^\ell\mid ia_i$} for all $i=1,\dots,d-1$, or equivalently, that $\ell v(p)\leqslant v(ia_i)$. Suppose that this is not the case, so there is some $i=1,\dots,d-1$ for which $\ell v(p)>v(ia_i)$. Then the Newton polygon of $f'$ will be the lower convex hull of the points $(j-1,v(ja_j))$ for $j=1,\dots,d-1$, together with the point $(d-1,v(d))$, so the Newton polygon will have a line segment of slope
    \[m=\max_{j=1,\dots,d-1}\frac{v(d)-v(ja_j)}{d-j}\geqslant \frac{v(d)-v(ia_i)}{d-i}=\frac{\ell v(p)-v(ia_i)}{d-i}>0.\]
Hence $f'$ will have a root $\ga\in\Kbar$ with $v(\ga)=-m<0$, which contradicts the fact that $|\ga|_v\leqslant 1$ for all critical points $\ga$ of $f$. Therefore no such $i$ can exist, so $f'(z)=p^\ell g(z)$ for some $g\in\Ocal[z]$, and thus $|f'(x)|_v<1$ for all $x\in\Dbar(0,1)$. Now, the fact that $f$ has integral coefficients and that $f'\!\left(\Dbar(0,1)\right)\subseteq D(0,1)$ shows that we can apply Lemma \ref{lem:periodic-disk-has-unique-attracting-periodic-point} to see that any open disk of unit radius that is contained in $\Dbar(0,1)$ and periodic under $f$ will contain a unique attracting periodic point $P$ (of the same period of the disk) that attracts every point of the disk. That is, if $D(c,1)\subseteq\Dbar(0,1)$ is periodic under $f$ with period $\ell\geqslant1$, then there is a unique point $P\in D(c,1)$ that is periodic under $f$ with period $\ell$ and for which 
    \[\lim_{n\to\infty}f^{\ell n}(x)=P\]
for any $x\in D(c,1)$. We now fix a critical point $\ga\in\Kbar$ of $f$ and let $U=D(\ga,1)$ be the residue class containing $\ga$. We split into various cases depending on whether $U$ is wandering or preperiodic under $f$. 

\underline{Case 1}: First suppose that $U$ is wandering. Then $a\in f^n(U)$ for at most one value of $n\geqslant0$, so there is some integer $N\geqslant0$ such that $a\notin f^n(U)=D(f^n(\ga),1)$ for all $n\geqslant N$. But this means that $|f^n(\ga)-a|_v\geqslant1$, so $|f^n(\ga)-a|_v=1$ because $\ga,a\in\Ocal$ and $f(\Ocal)=\Ocal$, and thus $v(f^n(\ga)-a)=0$ for all $n\geqslant N$. Hence the limit in equation (\ref{eq:desired-limit}) is $0$. 

\underline{Case 2:} We now assume that $U$ is not wandering, so it is preperiodic under $f$. Then there is some $m\geqslant0$ such that $V=f^m(U)=D(f^m(\ga),1)\subseteq\Dbar(0,1)$ is periodic under $f$ with period $\ell$, and thus our earlier observation shows that $V$ contains a unique attracting periodic point $P\in V$ of period $\ell$ that attracts all of $V$. We split Case 2 into two cases depending on whether $a$ is in the orbit of $P$ under $f$ or not.

\underline{Case 2a:} We first consider the case where $a$ is in the orbit of $P$ under $f$, so $a=f^j(P)$ for some $j=0,1,\dots,\ell-1$. Then $a$ is a periodic point of $f$ of period $\ell$, and we let $\la=\!\left(f^\ell\right)'\!(a)$ be the multiplier of $a$. Note that if
\begin{equation}
    0=\la=\!\left(f^\ell\right)'(a)=\prod_{i=0}^{\ell-1}f'(f^i(a))\label{eq:multiplier}
\end{equation}
then $f'(f^i(a))=0$ for some $i=0,\dots,\ell-1$. But if $f^i(a)=\ga$ for some $i>0$ and $\ga$ a critical point of $f$, then choosing $n>0$ such that $\ell n>i$ shows that $a=f^{n\ell}(a)=f^{n\ell-i}(\ga)$, and thus $a$ is in the postcritical set of $f$. Thus $\la=0$ is not possible, so we must have that $|\la|_v>0$. The fact that $|f'(x)|_v<1$ for all $x\in\Dbar(0,1)$ and equation (\ref{eq:multiplier}) then shows that $|\la|_v<1$ because $a\in\Dbar(0,1)$, so $0<|\la|_v<1$, i.e., $0<v(\la)<\infty$. Thus $a$ is an attracting but not superattracting periodic point of $f$. We now claim that if $x\in D(a,|\la|_v)$, then $\left|f^\ell(x)-a\right|_v=|\la|_v|x-a|_v$. We know that we can write $f^\ell(z)-a$ as
    \[f^\ell(z)-a=f^\ell(z)-f^\ell(a)=\la(z-a)+\sum_{i=2}^{d}b_i(z-a)^i,\]
where $b_2,\dots,b_{d-1}\in\Ocal$ by the comments at the beginning of the proof of Lemma \ref{lem:image-of-disc-under-monic-poly}, and where $b_d=1$ since $f$ is monic. If $x\in D(a,|\la|_v)$ is distinct from $a$ then 
    \[|\la|_v|x-a|_v>|x-a|_v^2\geqslant |b_i|_v|x-a|_v^i\]
for all $i=2,\dots,d$, so it follows that
    \begin{align*}
    \left|f^\ell(x)-a\right|_v&=\left|\la(x-a)+\sum_{i=2}^{d}b_i(x-a)^i\right|_v\\
    &=\max\!\left\{|\la|_v|x-a|_v,|b_2|_v|x-a|_v^2,\dots,|x-a|_v^d\right\}\\
    &=|\la|_v|x-a|_v.
    \end{align*} 
From this, we also see that if $i>0$ is an integer then
    \[\left|f^{i\ell}(x)-a\right|=|\la|_v^i|x-a|_v\]
for all $x\in D(a,|\la|_v)$. We now claim that there is some constant $c_0\in\Rbb$ such that for all sufficiently large $n$, we have that 
    \[v(f^n(\ga)-a)=\begin{cases}
        0, &\ell\nmid n-m-j,\\
        iv(\la)+c_0, \qquad&\dfrac{n-m-j}{\ell}=i\in\Zbb.
    \end{cases}\]
First consider the case where $\ell\nmid n-m-j$, so we can write $n-m-j=\ell s+t$ for some integer $s>0$ and $0<t<\ell$. In this case we have that
    \[f^n(\ga)=f^{\ell s+m+j+t}(\ga)\in f^{\ell s+m+j+t}(U)=f^{\ell s+j+t}(V)=f^{j+t}(V),\]
where we've used the fact that $V=f^m(U)$ is periodic under $f$ with period $\ell$. But 
    \[V=f^m(U)=f^m(D(\ga,1))=D(f^m(\ga),1),\]
so the fact that $P\in V$ shows that $V=D(P,1)$ and thus 
    \[f^n(\ga)\in f^{j+t}(D(P,1))=D(f^{j+t}(P),1)=D(f^t(a),1).\]
We now claim that $D(f^t(a),1)$ is disjoint from $D(a,1)$, which will shows that $|f^n(\ga)-a|_v=1$. Suppose that this is not the case, so $f^t(D(a,1))=D(a,1)$, and thus 
    \begin{multline*}
        f^{j+t}(V)=f^{j+t}(D(P,1)) =f^t(D(f^j(P),1)) =f^t(D(a,1)) =D(a,1)\\
        =D(f^j(P),1) =f^j(D(P,1)) =f^j(V).
    \end{multline*}
Hence 
    \[f^{m+j+t}(U)=f^{j+t}(V)=f^j(V)=f^{m+j}(U),\]
so the fact that $0<t<\ell$ shows that this contradicts the assumption that $U$ is preperiodic under $f$ with period $\ell$. Hence $D(f^t(a),1)\ni f^n(\ga)$ is disjoint from $D(a,1)$, so $|f^n(\ga)-a|_v=1$. We now deal with the case where $\ell|n-m-j$. Our previous observation that $V$ contains a unique periodic point $P$ of $f$ that attracts every element of $V$, and the fact that $f^m(\ga)\in V$, shows that
    \[\lim_{i\to\infty}f^{\ell i}(f^m(\ga))=P,\]
and thus
    \[\lim_{i\to\infty}f^{\ell i+m+j}(\ga)=a.\]
Hence there is some $i_0>0$ such that
    \[\left|f^{\ell i_0+m+j}(\ga)-a\right|_v<|\la|_v,\]
so $f^{\ell i_0+m+j}(\ga)\in D(a,|\la|_v)$. We know that $\ell i=n-m-j$ for some $i\in\Zbb$, and if we choose $n$ large enough so that $i>i_0$, then we have that
    \[\left|f^{n-m-j-\ell i_0}(x)-a\right|_v=\left|f^{\ell(i-i_0)}(x)-a\right|_v=|\la|_v^{i-i_0}|x-a|_v\]
for any $x\in D(a,|\la|_v)$. If we take $x=f^{\ell i_0+m+j}(\ga)\in D(a,|\la|_v)$, then this shows that
    \[\left|f^n(\ga)-a\right|_v=\left|f^{n-m-j-\ell i_0}\!\left(f^{\ell i_0+m+j}(\ga)\right)-a\right|_v=|\la|_v^{i-i_0}\left|f^{\ell i_0+m+j}(\ga)-a\right|_v,\]
so 
    \[v(f^n(\ga)-a)=(i-i_0)v(\la)+v\!\left(f^{\ell i_0+m+j}(\ga)-a\right)=iv(\la)+c_0,\]
where 
    \[c_0\coloneqq v\!\left(f^{\ell i_0+m+j}(\ga)-a\right)-i_0v(\la)<\infty\]
because $a$ is not in the postcritical set of $f$ and $0<v(\la)<\infty$. Thus if $n>0$ is sufficiently large then
    \[v(f^n(\ga)-a)=\begin{cases}
        0, &\ell\nmid n-m-j,\\
        iv(\la)+c_0, \qquad&\dfrac{n-m-j}{\ell}=i\in\Zbb,
    \end{cases}\]
which shows that $0\leqslant v(f^n(\ga)-a)\leqslant nv(\la)+c_0$ for all sufficiently large $n$, and thus
    \[0\leqslant\lim_{n\to\infty}\frac{v(f^n(\ga)-a)}{\psi(n)}\leqslant\lim_{n\to\infty}\frac{nv(\la)+c_0}{\psi(n)}=0.\]
It immediately follows that the desired limit is $0$. 

\underline{Case 2b:} We now consider the case where $a$ is not in the orbit of $P$ under $f$, or equivalently since $P$ is periodic under $f$ with period $\ell$, that $a\neq f^j(P)$ for all $j=0,\dots,\ell-1$. We know that
    \[\lim_{i\to\infty}f^{\ell i}(f^m(\ga))=P\]
because $P$ attracts every point of $V$ and $f^m(\ga)\in V$, so for any $j=0,\dots,\ell-1$ we have that
    \[\lim_{i\to\infty}f^{\ell i+m+j}(\ga)=f^j(P).\]
Hence for every $\ve>0$ there is some integer $I_\ve>0$ such that for all $j=0,\dots,\ell-1$, if $i\geqslant I_\ve$ then 
    \[f^{\ell i+m+j}(\ga)\in D(f^j(P),\ve).\]
The fact that $a\neq f^j(P)$ for all $j=0,\dots,\ell-1$ shows that there is some $0<r<1$ such that $a\notin D(f^j(P),r)$ for all $j=0,\dots,\ell-1$. This means that $f^{\ell i+m+j}(\ga)\notin D(a,r)$ for all $i\geqslant I_r$ and $j=0,\dots,\ell-1$, for if this wasn't true then $f^{\ell i+m+j}(\ga)$ would be in $D(f^j(P),r)\cap D(a,r)$ for some $i\geqslant I_r$ and $j=0,\dots,\ell-1$, which would show that $D(f^j(P),r)=D(a,r)$ and thus $a\in D(f^j(P),r)$. Hence if $j=0,\dots,\ell-1$ and $i\geqslant I_r$ then 
    \[r\leqslant\left|f^{\ell i+m+j}(\ga)-a\right|_v\leqslant1,\]
which means that if $n$ is sufficiently large, i.e., greater then $\ell I_r+m$, then $r\leqslant|f^n(\ga)-a|_v\leqslant 1$. It then follows that $0\leqslant v(f^n(\ga)-a)\leqslant -\log(r)$, so the desired limit in equation (\ref{eq:desired-limit}) is $0$. This deals with all the possible cases, so it follows that 
    \[\lim_{n\to\infty}\frac{v(f^n(\ga)-a)}{\psi(n)}=0.\qedhere\]
\end{proof}

\subsection{An Analysis Lemma}
\hspace{\parindent} The following lemma is used in the proof of Proposition \ref{prop:e_n-bound-in-composite-case}. The author believes this is a standard result, but due to lack of a source we include the proof. 

\begin{lemma}\label{lem:Analysis-lemma}
Let $(a_n)_{n\geqslant1}$ be a sequence of nonnegative real numbers and suppose that $b>1$ is a real number. Then
    \[\lim_{n\to\infty}\frac{a_n}{b^n}=0\quad\implies\quad\lim_{n\to\infty}\frac{1}{b^n}\sum_{k=1}^{n}a_k=0.\]
\end{lemma}

\begin{proof}
Fix $\ve>0$ and choose some $0<\de<\ve$. Because $a_n/b^n\to0$ as $n\to\infty$, we can find some integer $N_1\geqslant1$ such that if $n>N_1$ then 
    \[\frac{a_n}{b^n}<\frac{b-1}{b}\de.\]
Now choose an integer $N_2\geqslant1$ such that if $n>N_2$ then
    \[\frac{1}{b^n}\sum_{k=1}^{N_1}a_k<\ve-\de.\]
Let $N=\max\{N_1,N_2\}$. Then if $n>N$, we have that
\begin{align*}
    \frac{1}{b^n}\sum_{k=1}^{n}a_k
    &=\frac{1}{b^n}\sum_{k=1}^{N_1}a_k+\frac{1}{b^n}\sum_{k=N_1+1}^{n}a_k\\
    &=\frac{1}{b^n}\sum_{k=1}^{N_1}a_k+\sum_{k=N_1+1}^{n}\frac{a_k}{b^
    k}\frac{1}{b^{n-k}}\\
    &<\ve-\de+\sum_{k=N_1+1}^{n}\left(\frac{b-1}{b}\de\right)\frac{1}{b^{n-k}}\\
    &=\ve-\de+(b-1)\de\sum_{k=N_1+1}^{n}b^{k-1-n}\\
    &=\ve-\de+(b-1)\de\frac{1-b^{N_1-n}}{b-1}\\
    &=\ve-\de+\de\left(1-b^{N_1-n}\right)\\
    &<\ve-\de+\de=\ve.
\end{align*}
But $\ve$ was arbitrary, so this shows that
    \[\lim_{n\to\infty}\frac{1}{b^n}\sum_{k=1}^{n}a_k=0.\qedhere\]
\end{proof}

\subsection{Some Lemmas About Preimage Sets}\label{sect:some-lemmas-about-preimage-sets}
\hspace{\parindent} The following facts about preimage sets of functions are used in the proofs of some of the main theorems in this paper. The first one allows us to consider far fewer cases then initially seems necessary when looking at what possible directions elements in an iterated preimage set associated to $f$ can lie in.

\begin{lemma}\label{lem:preim-tree-either-has-size-1-or-is-infinite}
Suppose that $X$ is an infinite set and that $f:X\to X$ is a surjective function with the property that there is at most one point $x\in X$ for which $\#f^{-1}(x)=1$. Then for any $a\in X$, either $\#T_\infty(f,a)=1$, or $\#T_\infty(f,a)=\infty$. 
\end{lemma}
\begin{proof}
First assume that there is some $a'\in T_\infty(f,a)$ for which $\#f^{-1}(a')=1$. If $f^{-1}(a')=\{a'\}$ then $f(a')=a'$, which means that $a'=a$ since $a'=f^n(a')=a$ for some $n$ by definition. Thus $\#T_\infty(f,a)=1$. If $f^{-1}(a')=b'\neq a'$ then we replace $a$ with $a'$, since it will clearly suffice to show that $T_\infty(f,a')$ is infinite in order to show that $T_\infty(f,a)$ is. We can thus assume that there is some $b\in f^{-1}(a)$ distinct from $a$, and that for all $y\in T_\infty(f,a)$ distinct from $a$, $\#f^{-1}(y)\geqslant2$.

Now choose some $b\in f^{-1}(a)$ distinct from $a$. If $b\in f^{-1}(b)$ then $b\in f^{-1}(b)\cap f^{-1}(a)$, which is not possible since $b\neq a$. Hence $b\notin f^{-1}(b)$. Because $\#f^{-1}(b)\geqslant2$, it follows that we can choose some $\be_1\in f^{-1}(b)$ distinct from $a,b$. If $b\in f^{-1}(\be_1)$ then $b\in f^{-1}(\be_1)\cap f^{-1}(a)$, which is not possible since $\be_1\neq a$, and if $\be_1\in f^{-1}(\be_1)$ then $\be_1\in f^{-1}(b)\cap f^{-1}(\be_1)$, which is also not possible since $b\neq\be_1$. We now claim that for any $n\geqslant1$, there exists distinct $\be_1,\dots,\be_n\in T_\infty(f,a)\setminus\{a,b\}$ such that 
\begin{enumerate}[label=(\arabic*)]
    \item $\be_1\in f^{-1}(b)$ and $\be_i\in f^{-1}(\be_{i-1})$ for all $i=2,\dots,n$;
    \item and $b,\be_i\notin f^{-1}(\be_i)$ for all $i=1,\dots,n$.
\end{enumerate}
We know this is true for $n=1$, so suppose that it is true for some $n\geqslant1$. Because $b,\be_n\notin f^{-1}(\be_n)$, and $\#f^{-1}(\be_n)\geqslant2$, we know that we can find some $\be_{n+1}\in f^{-1}(\be_n)$ distinct from $a,b,\be_n$. If $\be_1=\be_{n+1}$ then $\be_1\in f^{-1}(b)\cap f^{-1}(\be_n)$, which is not possible since $b\neq\be_n$. If $\be_i=\be_{n+1}$ for any $i=2,\dots,n$ then $\be_i\in f^{-1}(\be_{i-1})\cap f^{-1}(\be_n)$, which is also not possible since $\be_{i-1}\neq\be_n$. Hence $\be_1,\dots,\be_{n+1}\in T_\infty(f,a)\setminus\{a,b\}$ are distinct. If $b\in f^{-1}(\be_{n+1})$ then $b\in f^{-1}(\be_{n+1})\cap f^{-1}(a)$, which is not possible since $a\neq\be_{n+1}$, and if $\be_{n+1}\in f^{-1}(\be_{n+1})$ then $\be_{n+1}\in f^{-1}(\be_{n})\cap f^{-1}(\be_{n+1})$, which is also not possible because $\be_n\neq\be_{n+1}$. Thus $\be_1,\dots,\be_{n+1}\in T_{\infty}(f,a)$ satisfy the necessary conditions. It follows that for any $n\geqslant1$, there exists distinct $\be_1,\dots,\be_n\in T_\infty(f,a)$, so $\#T_\infty(f,a)=\infty$.  
\end{proof}

The next lemma allows us to choose a subset of an iterated preimage set associated to $f$ that avoids certain elements. If $f:X\to X$ is a function and $x\in X$, then we set 
    \[f^\infty(x)=\{x,f(x),f^2(x),\dots\}.\] 

\begin{lemma}\label{lem:avoiding-elements-in-infinite-tree}
Let $X$ be an infinite set, let $f:X\to X$ be surjective, let $x\in X$, and suppose that $T_\infty(f, x)$ is an infinite set. Then the following are true:
\begin{enumerate}[label=(\alph*)]
    \item\label{avoiding-elements-in-infinite-tree-a} For any $y\in T_\infty(f,x)$, the set $f^\infty(y)\cap T_\infty(f,x)$ is finite.
    \item\label{avoiding-elements-in-infinite-tree-b} For any $y_1,\dots,y_N\in T_\infty(f, x)$, there exists some $x_0\in T_\infty(f, x)$ such that $y_1,\dots,y_N\notin T_\infty(f, x_0)$. 
\end{enumerate}
\end{lemma}
\begin{proof}
\ref{avoiding-elements-in-infinite-tree-a}: Note that if $f^n(y)\in T_\infty(f, x)$ for some $n$, then there is some $k\geqslant0$ such that $f^{n+k}(y)=x$, and thus if $0\leqslant m<n$ then 
    \[f^{n-m+k}(f^m(y))=f^{n+k}(y)=x.\]
Hence if $f^n(y)\in T_\infty(f, x)$, then $f^m(y)\in T_\infty(f, x)$ for all $0\leqslant m\leqslant n$. Thus if $f^\infty(y)\cap T_\infty(f, x)$ is infinite, then $f^n(y)\in T_\infty(f, x)$ for arbitrarily large $n$, which means that $f^\infty(y)\subseteq T_\infty(f, x)$. We now claim that $f^\infty(y)\subseteq T_\infty(f, x)$ implies that $f^\infty(y)$ is finite. Let $n\geqslant0$ be the smallest positive integer for which $f^n(y)=x$. Then $f^{n+1}(y)\in f^{\infty}(y)\subseteq T_\infty(f, x)$, so there is some $m\geqslant1$ such that 
    \[x=f^{n+m}(y)=f^m(f^n(y))=f^m(x).\]
Hence for any $k\geqslant 1$, $f^{n+mk}(y)=x$, so $f^{\infty}(y)=\{y,f(y),\dots,f^{n-1}(y),x,f(x),\dots,f^{m-1}(x)\}$ is finite, and thus $T_\infty(f, x)\cap f^{\infty}(y)$ cannot be infinite. 

\ref{avoiding-elements-in-infinite-tree-b}: Part \ref{avoiding-elements-in-infinite-tree-a} shows that for any $y_1,\dots,y_N\in T_\infty(f, x)$, we have that 
    \[T_\infty(f, x)\cap\bigcup_{i=1}^{N}f^{\infty}(y_i)=\bigcup_{i=1}^{N}T_\infty(f, x)\cap f^{\infty}(y_i)\]
is finite. Since $T_\infty(f, x)$ is infinite, it follows that we can choose some $x_0\in T_\infty(f, x)$ that is not in $f^{\infty}(y_i)$ for all $i=1,\dots,N$. Hence $y_1,\dots,y_N\notin T_\infty(f, x_0)\subseteq T_\infty(f, x)$, since if $y_i\in T_\infty(f, x_0)$ then $x_0\in f^{\infty}(y_i)$ by definition. 
\end{proof}

\nocite{*}
\printbibliography

\end{document}